\documentclass[runningheads]{llncs}

\usepackage{amsmath}
\usepackage{amssymb}
\usepackage{xspace}
\usepackage{xcolor}
\definecolor{darkgreen}{rgb}{0,0.5,0}
\colorlet{COLOR1}{orange}
\colorlet{COLOR2}{darkgreen}
\colorlet{COLOR3}{blue}
\definecolor{COLOR4}{rgb}{0.5,0,0.5}
\colorlet{COLOR5}{red}
\colorlet{COLOR6}{black}
\usepackage{jlcode}
\usepackage{listings}
\usepackage[labelformat=simple]{subcaption}

\usepackage{csvsimple}
\usepackage{booktabs}
\usepackage{siunitx}
\usepackage{pgfplots}
\usetikzlibrary{spy}
\usepackage{hyperref}

\newcommand{\mc}[1]{\multicolumn{1}{c}{#1}}
\newcommand{\R}{\ensuremath{\mathbb{R}}\xspace}
\newcommand{\X}{\ensuremath{\mathcal{X}}\xspace}
\newcommand{\U}{\ensuremath{\mathcal{U}}\xspace}
\newcommand{\B}{\ensuremath{\mathcal{B}}\xspace}
\newcommand{\CH}{\ensuremath{\mathit{CH}}\xspace}
\newcommand{\Y}{\ensuremath{\mathcal{Y}}\xspace}

\newcommand{\V}{\ensuremath{\mathcal{V}}\xspace}
\newcommand{\reach}{\ensuremath{\mathcal{R}}\xspace}

\let\deltatmp\delta  

\begin{document}

\title{Conservative Time Discretization: \\ A Comparative Study}

\author{%
Marcelo Forets\inst{1}\orcidID{0000-0002-9831-7801}
\and \\
Christian Schilling\inst{2}\orcidID{0000-0003-3658-1065}
}

\authorrunning{M.\ Forets \and C.\ Schilling}

\institute{%
DMA, CURE, Universidad de la Rep\'ublica, Uruguay
\\
\email{mforets@gmail.com} \\
\and
Aalborg University, Aalborg, Denmark
\\
\email{christianms@cs.aau.dk}
}

\maketitle

\begin{abstract}
We present the first review of methods to overapproximate the set of reachable states of linear time-invariant systems subject to uncertain initial states and input signals \emph{for short time horizons}.
These methods are fundamental to state-of-the-art reachability algorithms for long time horizons, which proceed in two steps:
First they use such a method to discretize the system for a short time horizon, and then they efficiently obtain a solution of the new discrete system for the long time horizon.
Traditionally, both qualitative and quantitative comparison between different reachability algorithms has only considered the combination of both steps.
In this paper we study the first step in isolation.
We perform a variety of numerical experiments for six fundamental discretization methods from the literature.
As we show, these methods have different trade-offs regarding accuracy and computational cost and, depending on the characteristics of the system, some methods may be preferred over others.
We also discuss preprocessing steps to improve the results and efficient implementation strategies.

\end{abstract}

\section{Introduction} \label{sec:introduction}

\begin{figure}[tb!]
	\centering
	\begin{subfigure}[b]{0.49\textwidth}
		\centering
		\includegraphics[width=\linewidth,keepaspectratio]{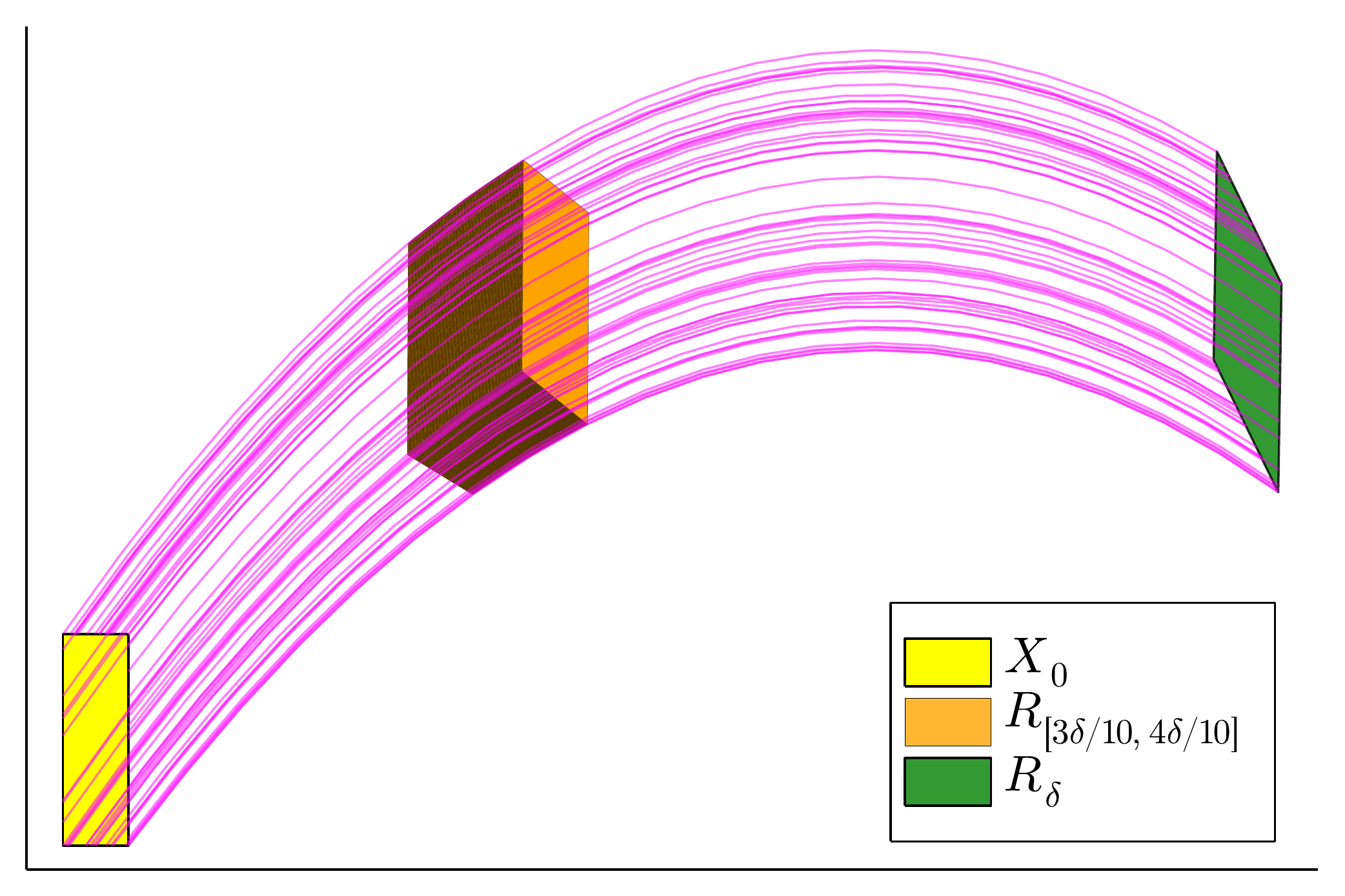}
		\caption{Sketch of the reachability problem.}
		\label{fig:intro_exact}
	\end{subfigure}
	\begin{subfigure}[b]{0.49\textwidth}
		\centering
		\includegraphics[width=\linewidth,keepaspectratio]{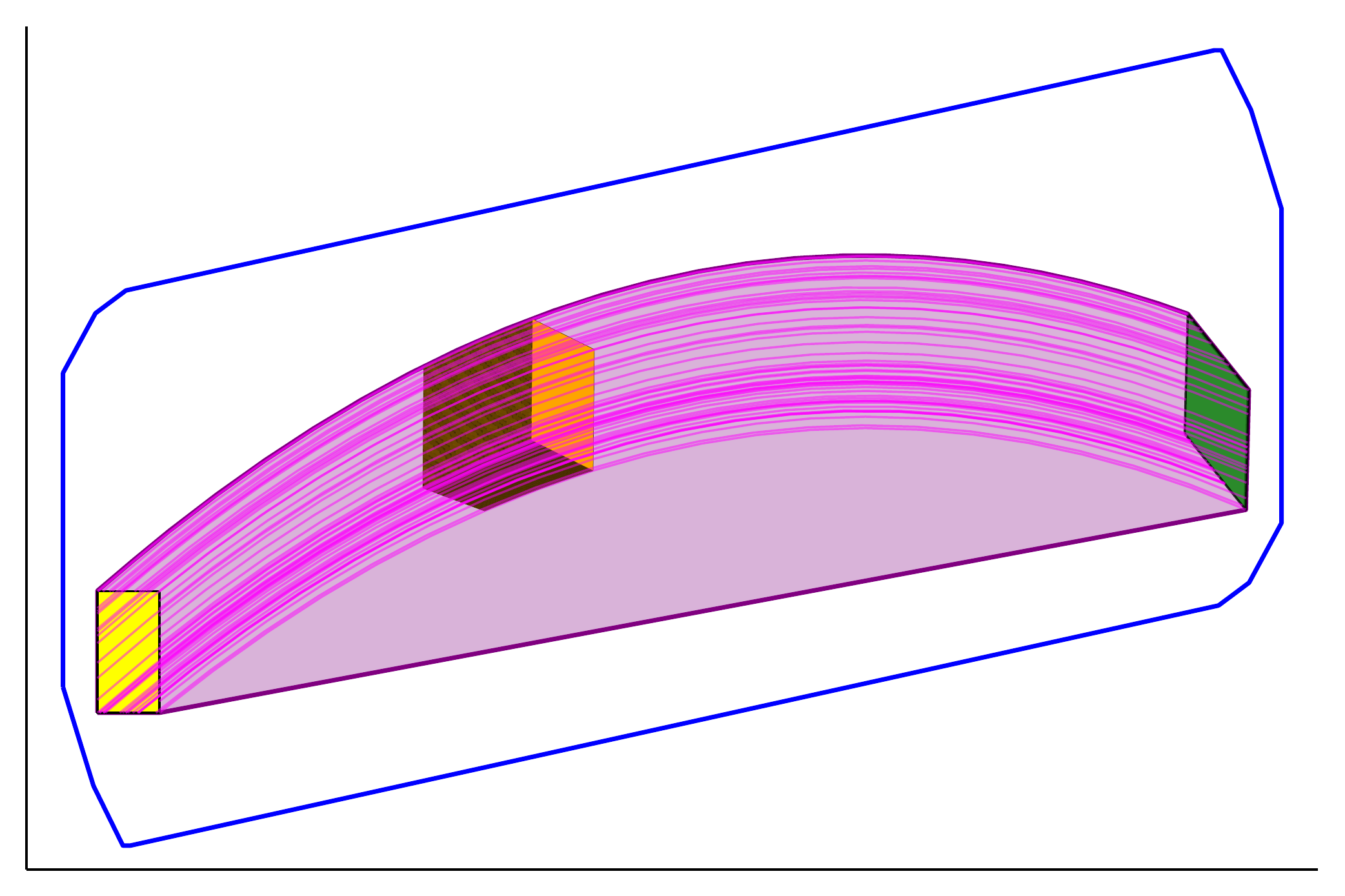}
		\caption{Two possible solutions for $\Omega_0$.}
		\label{fig:intro_approximate}
	\end{subfigure}
	
	\begin{subfigure}[b]{\textwidth}
		\centering
		\includegraphics[width=\linewidth,height=5cm,keepaspectratio]{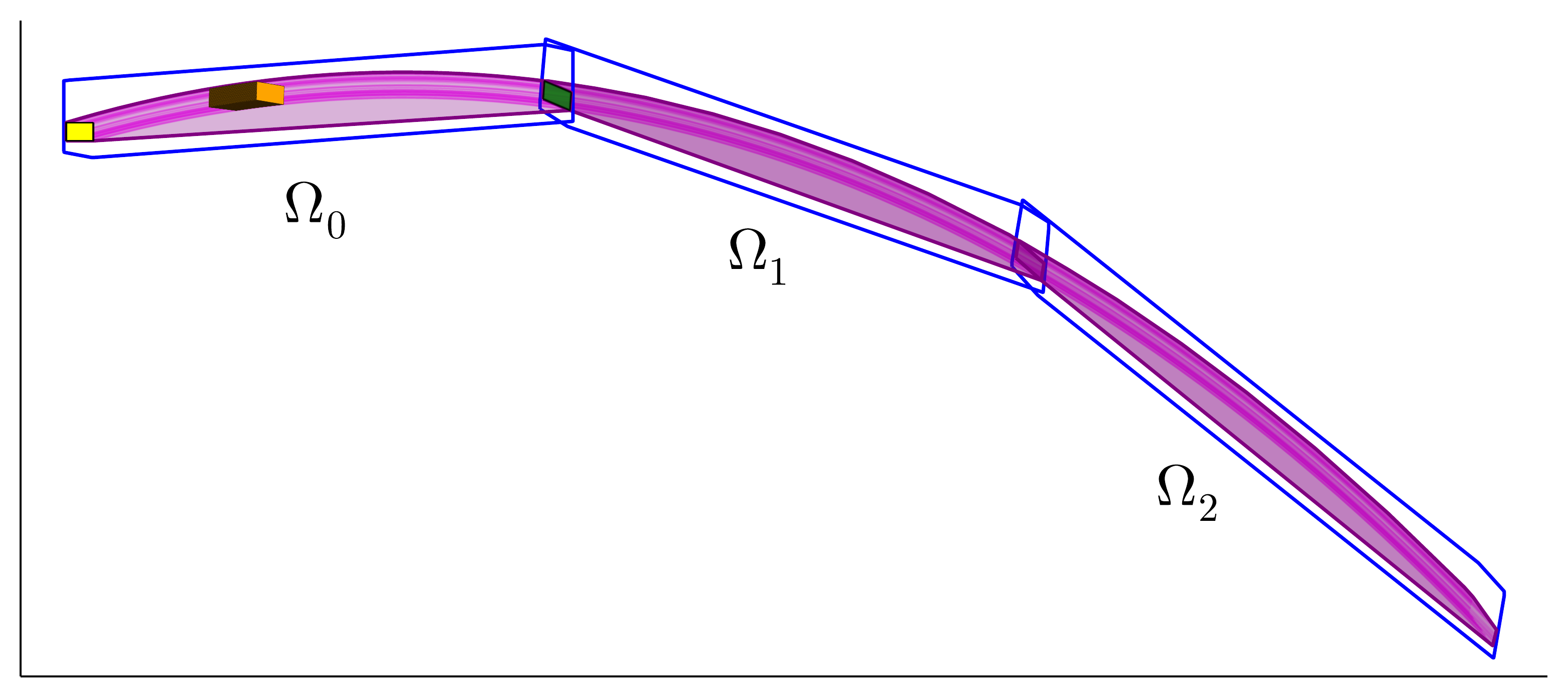}
		\caption{First three steps of the iteration in a reachability algorithm.}
		\label{fig:intro_recurrence}
	\end{subfigure}
	\caption{\textbf{Left:} Starting from the initial states $\X_0$ (yellow), we want to enclose the states reachable within $[0, \delta]$, $\reach_{[0, \delta]}$, covering the random trajectories (magenta).
	\textbf{Right:} The purple set is the smallest convex solution for $\Omega_0$.
	The blue set is a zonotope solution for $\Omega_0$ computed with the method in \cite{AlthoffSB07}.
	\textbf{Bottom:} The first three sets of the discretized system computed by a reachability algorithm.
	}
\end{figure}

We study the fundamental problem of reachability for a system of linear differential equations.
Given a set of initial states $\X_0 \subseteq \R^n$, we are interested in the set of states that can be reached by any trajectory up to some time horizon.

The classical analysis approach is numerical simulation, which has several drawbacks.
First, one can only simulate finitely many trajectories, but a system has infinitely many trajectories.
Second, a simulated trajectory is only available at finitely many discrete points in time.
Third, even for these discrete points in time, a standard simulation is not guaranteed to be exact.
For all these reasons, simulation may miss critical behaviors, which can lead to wrong conclusions.

In contrast, reachability algorithms construct a finite sequence of \emph{sets} covering \emph{all} possible \emph{continuous} trajectories of the system.
In general it is not possible to represent the exact set of reachable states, but for linear systems one can obtain arbitrary-precision approximations as a union of convex sets \cite{le2010reachability,AlthoffFG20}.

The fundamental procedure behind all modern reachability algorithms \cite{ARCH-COMP} consists of two stages.
In the first stage, the system is discretized, i.e., approximated by a discrete system.
This requires to find a Euclidean set $\Omega_0 \subseteq \R^n$ that includes all states reachable starting from any initial state $x(0) \in \X_0$ up to a small time horizon $\delta$, as illustrated in Fig.~\ref{fig:intro_exact}.
Two possible solutions for $\Omega_0$ are presented in Fig.~\ref{fig:intro_approximate}.
In the second stage, the set $\Omega_0$ is propagated forward in time until the time horizon is reached, which is sketched in Fig.~\ref{fig:intro_recurrence}.
Thus computing a precise $\Omega_0$ is important for the precision of the overall result.

All state-of-the-art reachability tools such as CORA \cite{Althoff15}, the continuous version of Hylaa \cite{BakTJ19}, HyPro \cite{SchuppAMK17}, JuliaReach \cite{BogomolovFFPS19}, and SpaceEx \cite{FrehseGDCRLRGDM11} follow this procedure of discretizing and computing an approximation $\Omega_0$.
Over the years, multiple methods to obtain such sets $\Omega_0$ have been proposed  \cite{ChutinanK99,AsarinDMB00,Girard05,AlthoffSB07,le2010reachability,FrehseGDCRLRGDM11,BogomolovFFVPS18,althoff2020reachability}.
These methods are tailored toward different set representations as a requirement of the second stage of the reachability algorithm.
The two prevalent options are a concrete representation with a zonotope and a lazy representation based on the support function.
Nevertheless, the different sets $\Omega_0$ solve the same problem and can be partially interchanged between different algorithms.
However, to our knowledge, these methods have never been compared to each other.

In this article we study different methods to compute $\Omega_0$.
We introduce this problem formally in the next section.
Our study has these particular goals:
\begin{itemize}
	\item Present the methods in a unified way (Section~\ref{sec:methods}).
	
	\item Discuss potential gains of system transformations (Section~\ref{sec:transformation}).
	
	\item Discuss aspects of an efficient implementation (Section~\ref{sec:implementation}).
	
	\item Assess the effect of different system characteristics and evaluate the methods on different systems (Section~\ref{sec:evaluation}).
\end{itemize}

\section{Problem statement}\label{sec:problem}

In this work we study linear time-invariant (LTI) systems, which have the state-space form
\begin{equation}\label{eq:LTI_general}
	\dot{x}(t) = A x(t) + u(t),
\end{equation}
where $x$ is an $n$-dimensional state vector, $A \in \R^{n \times n}$ is the flow matrix, and $u$ is a bounded but arbitrarily varying input signal that belongs to the input domain $\U \subseteq \R^n$, i.e., $u(t) \in \U$ for all $t \geq 0$.
It is common to associate a linear map $B$ with $u$ in \eqref{eq:LTI_general}, but this map can be absorbed in $\U$ (via $\U \mapsto B \U$, which is no restriction in practice because common set representations for \U such as zonotopes are closed under linear maps).
An LTI system is homogeneous if $\U = \{0\}$ and heterogeneous otherwise.
We consider initial-value problems for system~\eqref{eq:LTI_general} where the initial state $x(0) \in \mathbb{R}^n$ is taken from a set of initial states $\X_0 \subseteq \R^n$.
The solution of \eqref{eq:LTI_general} for a particular initial state $x_0$ and input signal $u$ is the trajectory $\xi_{x_0, u}(t)$, which is a function of time.
The set of solutions at time $t$ is the set of reachable states
\begin{equation}
	\reach_t = \{ \xi_{x_0, u}(t) : x_0 \in \X_0, u \in \U \},
\end{equation}
and we generalize this set to time intervals
\begin{equation}
	\reach_{[t_0, t_1]} = \{ \xi_{x_0, u}(t) : x_0 \in \X_0, u \in \U, t \in [t_0, t_1] \}.
\end{equation}

The time-bounded reachability problem asks to compute the set of reachable states $\reach_{[0, T]}$ up to a time horizon $T$.
Solving this problem exactly for LTI systems is generally not possible \cite{LafferrierePY01}.
Hence an overapproximation $\Omega \supseteq \reach_{[0, T]}$ is sought in practice.
The common approach to compute the set $\Omega$ proceeds in two stages.

In the first stage, conservative time discretization is applied to compute a set $\Omega_0$ with the property that it encloses the exact reachable states for an initial time interval $[0, \delta]$, i.e., $\Omega_0 \supseteq \reach_{[0, \delta]}$; usually one chooses a small value for $\delta$, much smaller than $T$.
We illustrate two possible choices for $\Omega_0$ in Fig.~\ref{fig:intro_approximate}.
For efficiency reasons, in practice one restricts $\Omega_0$ to convex sets.
Under this restriction, the smallest set $\Omega_0$ is the convex hull of $\reach_{[0, \delta]}$ (purple set in Fig.~\ref{fig:intro_approximate}).

\smallskip

In a similar fashion to $\Omega_0$, the second stage also requires a discretized input set $\V$, which encloses the trajectories of system~\eqref{eq:LTI_general} but starting from the origin ($x(0) = 0$), up to the time horizon $\delta$.
Since the set $\U$ is often assumed to have a simple shape, the computation of $\V$ is straightforward and mostly identical across the different discretization methods, and we do not discuss it further.

The second stage propagates the set $\Omega_0$ through consecutive time intervals until reaching the time horizon $T$, which we sketch briefly because the technicalities are not of interest.
The sequence of sets $\Omega_k$, $k \geq 0$, is given by the following recurrence, where $\oplus$ denotes the Minkowski sum, $\X \oplus \Y = \{x + y : x \in \X, y \in \Y\}$:
\begin{equation}
	\Omega_k = e^{A \delta k} \Omega_0 \oplus \bigoplus_{j=0}^{k-1} e^{A \delta j} \V
\end{equation}

We finally define our approximation of the reachable states (see Fig.~\ref{fig:intro_recurrence}) as
\begin{equation}
	\Omega = \bigcup_{k=0} \Omega_k.
\end{equation}

We note that state-of-the-art reachability algorithms for LTI systems are wrapping-free, i.e., they do not accumulate errors over time; as such, the precision of these algorithms is mainly determined by the precision of $\Omega_0$ and $\V$.

The reader may wonder why reducing the problem of computing $\Omega \supseteq \reach_{[0, T]}$ to the problem of computing $\Omega_0 \supseteq \reach_{[0, \delta]}$, which is structurally equivalent, makes sense.
It is important to remark that there are good convex approximation methods for small enough values of $\delta$, but these methods fail for large values of $\delta$ (as we shall see later).
While $\Omega_0$ is convex, $\Omega$ is a union of convex sets (and thus typically not convex).
We also note the analogy to numerical simulation methods, which may also require that the time step is small enough, for instance by the well-known Courant-Friedrichs-Lewy condition.

In summary, the fundamental problem that we study here is:

\begin{problem}\label{prob:discretize}
	Given an LTI system and a time horizon $\delta$, find a set $\Omega_0 \supseteq \reach_{[0, \delta]}$.
\end{problem}

\section{Discretization methods}\label{sec:methods}

In this section, after fixing some common notation and mentioning approaches for nonlinear systems, we introduce various methods to obtain a conservative discretization $\Omega_0$ for solving Problem~\ref{prob:discretize}.
The presentation follows the chronological order in which these methods have been proposed.
The methods for computing $\Omega_0$ we consider here assume that $\X_0$ and $\U$ are at least compact and convex.
Under this condition it is known that $\reach_t$ is also compact and convex for any $t$ (unlike the sets $\reach_{[t_1, t_2]}$).
We end the section with pointers to related works outside the scope of this study.

\subsection{Notation}

In this article, when the norm of an operator is not specified, we assume a $p$-norm with $1 \leq p \leq \infty$.
We denote the ball of radius $\varepsilon$ in some $p$-norm and centered in the origin by $\B_\varepsilon$, and $\B_\varepsilon^p$ when the $p$-norm is relevant.
Let $\X, \Y \subseteq \R^n$ be sets.
The norm of \X is defined as $\Vert \X \Vert = \sup_{x \in \X} \Vert x \Vert$.
The symmetric interval hull of $\X$, written $\boxdot(\X)$, is the smallest hyperrectangle that contains $\X$ and is centrally symmetric in the origin.
We write $\CH(\X, \Y)$ to denote the convex hull of the union of $\X$ and $\Y$, and $\rho(d, \X)$ to denote the support function of $\X$ along direction $d \in \R^{n}$ (see for instance \cite{LazySets} for a formal definition and examples).
Given system \eqref{eq:LTI_general}, we define the state-transition matrix $\Phi = e^{A \delta}$ as the matrix exponential of $A \delta$. The following matrix function is also relevant:
\begin{equation} \label{eq:Phi2}
	\Phi_2(A, \delta) = \sum_{i=0}^\infty \frac{\delta^{i+2}}{(i+2)!} A^i.
\end{equation}
If $A$ is invertible, $\Phi_2(A, \delta) = A^{-2}(\Phi - I - A \delta)$. See Section~\ref{sec:matrix_functions} for the computation.

\subsection{Methods for nonlinear systems}\label{sec:methods_nonlinear}

In principle, we can apply reachability algorithms for nonlinear systems to LTI systems as a special case.
However, the methods developed specifically for LTI systems are much more precise and scalable.
Hence we do not include nonlinear approaches in our study and only briefly mention some relevant works below.

A polyhedral enclosure $\Omega_0$ can be computed for any dynamical system by choosing a set of normal directions $(d_i)_i$ and solving the corresponding optimization problems $\max_{t \in [0, \delta]} d_i^T x(t)$ \cite{ChutinanK99}.
This scheme relies on a sound optimizer and the run time depends on the number of directions and is difficult to predict.
For LTI systems, the analytic expression of $x(t)$ can be used \cite[Section 5.3]{le2009reachability}.

The work in \cite{AsarinDG03}, describes how to deal with nondeterministic inputs for Lipschitz-continuous systems.
If $L$ is the Lipschitz constant, the reachable states of the heterogeneous system are contained in the bloated reachable states of the homogeneous system: $f_\mathit{heterog} \subseteq f_\mathit{homog} \oplus \B_\varepsilon$, where $\varepsilon = \frac{\Vert \U \Vert}{L}(e^{L \delta} - 1)$.

\subsection{Common structure of methods for linear systems}\label{sec:structure}

Examining the different methods that we study in the remainder of this section, the following common structure emerges:
\begin{equation} \label{eq:Omega0_general}
	\Omega_0 = \CH(\X_0, \Phi \X_0 \oplus \mathcal{H}) \oplus \mathcal{J},
\end{equation}
for suitable sets $\mathcal{H}$ and $\mathcal{J}$ (possibly empty).
In short, the idea is to compute the convex hull between the reachable states at time $0$ and at time $\delta$, $\X_0$ and $\Phi \X_0$.
That would suffice if the trajectories were just following straight lines.
To correct for the curvature, the bloating terms $\mathcal{H}$ and $\mathcal{J}$ need to be added.

\subsection{First-order d/dt method} \label{ref:firstorder_ddt}

The earliest work specifically designed for linear systems we are aware of was developed for the tool d/dt \cite{AsarinDMB00}.
This work only considers homogeneous systems.
The definition is
\begin{equation}\label{eq:method_homogeneous}
	\Omega_0 = \CH(\X_0, \Phi \X_0) \oplus \B_\varepsilon,
\end{equation}
where
\begin{equation}
	\varepsilon = \left( e^{\Vert A \Vert \delta} - 1 - \Vert A \Vert \delta\right)\Vert \X_0 \Vert - \frac{3}{8}\Vert A\Vert^2 \delta^2 \Vert \X_0 \Vert.
\end{equation}

We note that the authors also claim that their method can be used to obtain an underapproximation, but this is not correct (see Appendix~\ref{sec:underapproximation} for details).

\subsection{First-order zonotope method}\label{sec:zonotope}

The work in \cite{Girard05} computes a zonotope enclosure.
The idea is to cover $\CH(\X_0, \Phi \X_0)$ with a zonotope and then bloat with another zonotope (here: a ball in the infinity norm):
\begin{equation}
	\Omega_0 = \mathit{zonotope}(\CH(\X_0, \Phi \X_0)) \oplus \B_{\varepsilon}^\infty
\end{equation}
where
\begin{equation}
	\varepsilon = \left( e^{\Vert A \Vert_\infty \delta} - 1 - \Vert A \Vert_\infty \delta\right)\left( \Vert \X_0 \Vert_\infty + \frac{\Vert \U \Vert_\infty}{\Vert A \Vert_\infty} \right) + \delta \Vert \U \Vert_\infty.
\end{equation}

\subsection{Correction-hull method} \label{sec:correction_hull}

The method in \cite{AlthoffSB07} is designed for interval dynamics matrix $A$, which represents uncertain parameters and has the scalar matrix as a special case.
The resulting set, a zonotope, is constructed from interval linear maps, which is described in \cite[Theorem 4]{AlthoffSB07}.
The approach is based on truncating the Taylor series at a chosen order $p$, and this order must satisfy the following inequality:
\begin{equation}\label{eq:correction_hull_convergence}
	\alpha = \frac{\Vert A \Vert_\infty \delta}{p + 2} < 1.
\end{equation}

The method assumes that the input domain \U contains the origin.
If this is not the case, a simple transformation can bring the system to this form \cite[Section~3.2.2]{Althoff10}, which we describe in Section~\ref{sec:homogenize}.
Then we have the following definitions:
\begin{align}
	\Omega_0 &= \CH(\X_0, \Phi \X_0) \oplus F_p \X_0 \oplus G_p \U,
	\intertext{where the correction matrices $F$ and $G$ are}
	F_p &= E + \sum_{i=2}^p [\delta^i (i^{\frac{-i}{i-1}} - i^{\frac{-1}{i-1}}), 0] \frac{A^i}{i!} \\
	\intertext{and}
	G_p &= E \delta + \sum_{i=0}^p \frac{A^i \delta^{i+1}}{(i+1)!}.
	\intertext{The remainder matrix $E$ to bound $\sum_{i=p+1}^{\infty} A^i \delta^i / i!$ is a diagonal interval matrix based on results in \cite{Liou66}; here the assumption \eqref{eq:correction_hull_convergence} is relevant to make the geometric series $\sum_i \alpha^i$ converge. Define}
	E &= [-\varepsilon, \varepsilon] \mathbf{1}
	\intertext{where $\mathbf{1}$ is the $n \times n$ matrix filled with ones and}
	\varepsilon &= \frac{(\Vert A \Vert_\infty \delta)^{p+1}}{(p+1)!} \frac{1}{1 - \alpha}.
\end{align}

The method was later extended in \cite{AlthoffKS11} and in our implementation we use the remainder term $E$ from that work instead.

\subsection{First-order method}

The method in \cite{le2010reachability} uses a first-order approximation similar to \cite{Girard05} in Section~\ref{sec:zonotope}, but in contrast it is not restricted to zonotopes and the infinity norm:
\begin{equation}
	\Omega_0 = \CH(\X_0, \Phi \X_0 \oplus \delta \U \oplus \B_{\varepsilon})
\end{equation}
where
\begin{equation}
	\varepsilon = \left( e^{\Vert A \Vert \delta} - 1 - \Vert A \Vert \delta\right)\left( \Vert \X_0 \Vert + \frac{\Vert \U \Vert}{\Vert A \Vert} \right).
\end{equation}

\subsection{Forward-backward method}

The approach in \cite{FrehseGDCRLRGDM11} describes an optimization procedure similar in spirit to the one discussed in Section~\ref{sec:methods_nonlinear} but specialized for LTI systems. Here one only needs to optimize over a quadratic function.
First we define some auxiliary terms.
\begin{align}
	E_+ &= \boxdot(\Phi_2(|A|, \delta) \boxdot(A^2 \X_0)) \label{eq:Eplus} \\
	E_- &= \boxdot(\Phi_2(|A|, \delta) \boxdot(A^2 \Phi \X_0)) \label{eq:Eminus} \\
	E_\psi &= \boxdot(\Phi_2(|A|, \delta) \boxdot(A \U))  \label{eq:Epsi} \\
	\Y_\lambda &= (1 - \lambda) \X_0 \oplus \lambda \Phi \X_0 \oplus \lambda \delta \U \oplus (\lambda E_+ \cap (1 - \lambda) E_-) \oplus \lambda^2 E_\psi \nonumber
\end{align}

Here the term $E_+$ goes forward from $\X_0$ and the term $E_-$ goes backward from $\Phi \X_0$, and it is sufficient to consider their intersection.
The solution is then obtained by optimizing $\Y_\lambda$ (where the objective function is piecewise-linear for homogeneous systems and piecewise-quadratic for heterogeneous systems):
\begin{equation}\label{eq:method_forwardbackward}
	\Omega_0 = \CH(\bigcup_{\lambda \in [0, 1]} \Y_\lambda ).
\end{equation}

\subsection{Forward-only method} \label{sec:forward_only}

The work in \cite{BogomolovFFVPS18} uses a simplified version of \eqref{eq:method_forwardbackward} with only a forward approximation, which works without an optimization procedure.
It can be seen that
\begin{equation}
	\CH(\bigcup_{\lambda \in [0, 1]} \Y_\lambda) \subseteq \CH(\Y_0, \Y_1 \oplus E_+)
\end{equation}
and this method accordingly uses
\begin{equation}
	\Omega_0 = \CH(\X_0, \Phi \X_0 \oplus \delta \U \oplus E_\psi \oplus E_+),
\end{equation}
where $E_+$ and $E_\psi$ are defined in \eqref{eq:Eplus} and \eqref{eq:Epsi}.
Analogously one can define a ``backward-only'' method by using $E_-$ instead of $E_+$.

\subsection{Combining methods}

It should be noted that if $\Omega_0^{a}$ and $\Omega_0^{b}$ are two solutions to Problem~\ref{prob:discretize}, then their intersection $\Omega_0^{a} \cap \Omega_0^{b}$ is also a solution.
(The dual statement holds for underapproximations and their unions.)
It is hence possible to combine the different methods outlined above.
This idea was used in \cite{ForetsCP21} where the authors combined the ``forward-only'' method with the mentioned ``backward-only'' method; this method yields solutions that are closer to the ``forward-backward'' method than these methods individually but is more efficient than the latter.

\subsection{Application to high-dimensional systems}\label{sec:krylov}

The approach in \cite{althoff2020reachability} shows how to efficiently work with high-dimensional systems.
The construction of $\Omega_0$ is similar to the ``correction-hull'' method, including the focus on zonotopes as set representation.
The difference is that the structure of the solution is rewritten to rely on matrices in the Krylov subspace.
The idea is to compute two matrices $W$ and $H$ to approximate the effect of a vector $v$ on the matrix exponential $e^A$ without computing it:
\begin{equation}
	e^A v \approx \Vert v \Vert W e^H e_1.
\end{equation}

The matrix exponential $e^H$ can be computed efficiently, and there exist estimates to bound the above approximation error \cite{althoff2020reachability}.

In \cite{ForetsCP21} the authors demonstrated that the ``forward-only'' method can also be efficiently implemented with Krylov techniques.

\subsection{Application to time-varying systems}

Linear systems whose dynamics are time-varying due to uncertain parameters can be represented with interval matrices.
This setting differs from the one in \cite{AlthoffSB07} from Section~\ref{sec:correction_hull} where the system dynamics are uncertain but time-invariant.
Methods based on zonotopes to handle such systems are presented in \cite{Althoff10,AlthoffGK11}.

\section{Problem transformations}\label{sec:transformation}

In this section we shortly explain possible transformations of system \eqref{eq:LTI_general} to a normal form.
These transformations are simple but yield interesting results.

\smallskip

For illustration, we use a simple harmonic oscillator with inputs $f$,
\begin{equation}\label{eq:sdof}
	\ddot{y}(t) + \omega^2 y(t) = f,
\end{equation}
where $\omega^2 = 4\pi$, $m = 1$ and $y(t)$ is the unknown. This problem can be associated with a spring-mass system, where $y(t)$ is the elongation of the spring at time $t$ and $\omega$ is the natural frequency.
Bringing Eq.~\eqref{eq:sdof} to the first-order form of Eq.~\eqref{eq:LTI_general} with the change of variables $x(t) = [y(t), \dot{y}(t)]^T$, we obtain
\begin{equation} \label{eq:sdof_Amatrix}
	\dot{x}(t) = \begin{pmatrix}
		0 & 1 \\
		-4\pi  & 0
	\end{pmatrix} x(t) + \binom{0}{f}.
\end{equation}

\subsection{Homogenization}\label{sec:homogenize}

The first transformation, which was used in \cite[Section~3.2.2]{Althoff10} for the ``correction-hull'' method (although with a different goal) and in \cite{ForetsCP21}, expresses some of the inputs' effect with a fresh state variable.
For deterministic systems (i.e., where the input domain \U is a singleton), this allows to completely eliminate the inputs.
For proper nondeterministic systems, this extracts the ``central'' effect of the input signals into a state variable and only leaves the deviation from this central effect, effectively re-centering the input domain \U in the origin.
The transformation is motivated because \U is treated rather pessimistically in the methods to compute $\Omega_0$ (e.g., it may appear in the form of its norm $\Vert \U \Vert$).

We illustrate the issue in one dimension.
For the $p$-norms considered here and any $a \in \R$ we have $\Vert a \Vert = \Vert -a \Vert$.
Hence the norm of the singleton $\U_1 = \{ a \}$ is equivalent to the norm of the proper interval $\U_2 = [-a, a]$.
Thus $\Omega_0$ is identical in both cases, but since it must cover the $\U_2$ case, with many more possible behaviors, it is coarse for $\U_1$ and large $a$.

Now suppose that the input domain \U is centrally symmetric but not centered in the origin.
The idea of the transformation is to shift \U to the origin and add another state variable to account for this shift.
Formally, assume a heterogeneous system \eqref{eq:LTI_general}, an initial set $\X_0$, and a domain $\U$ centered in a point $c \neq 0$.
We define a new system $\dot{y}(t) = C y(t) + u(t)$, where we have the block matrix
\begin{equation}
	C = \begin{pmatrix} A & \frac{b}{\alpha} \\ 0 & 0\end{pmatrix}
\end{equation}
for some value $\alpha \neq 0$, $\Y_0 = \X_0 \times \{\alpha\}$ is the new initial set, and $\U \oplus \{-c\}$ is the new input domain.
The parameter $\alpha$ can be used to trade off the impact in the norm of $C$ (for methods where this term appears) and $\Y_0$.
This transformation increases the state dimension $n$ by one and removes the input dimension.
The modification of the initial states (i.e., the construction of $\Y_0$ from $\X_0$) is efficient for typical set representations of $\X_0$.
Likewise, projecting away the auxiliary dimension in the end is efficient for common set representations since this dimension is independent of all other dimensions and has the constant value $1$.

\begin{figure}[t]
	\centering
	\begin{subfigure}[b]{0.49\textwidth}
		\centering
		\includegraphics[width=\linewidth,keepaspectratio]{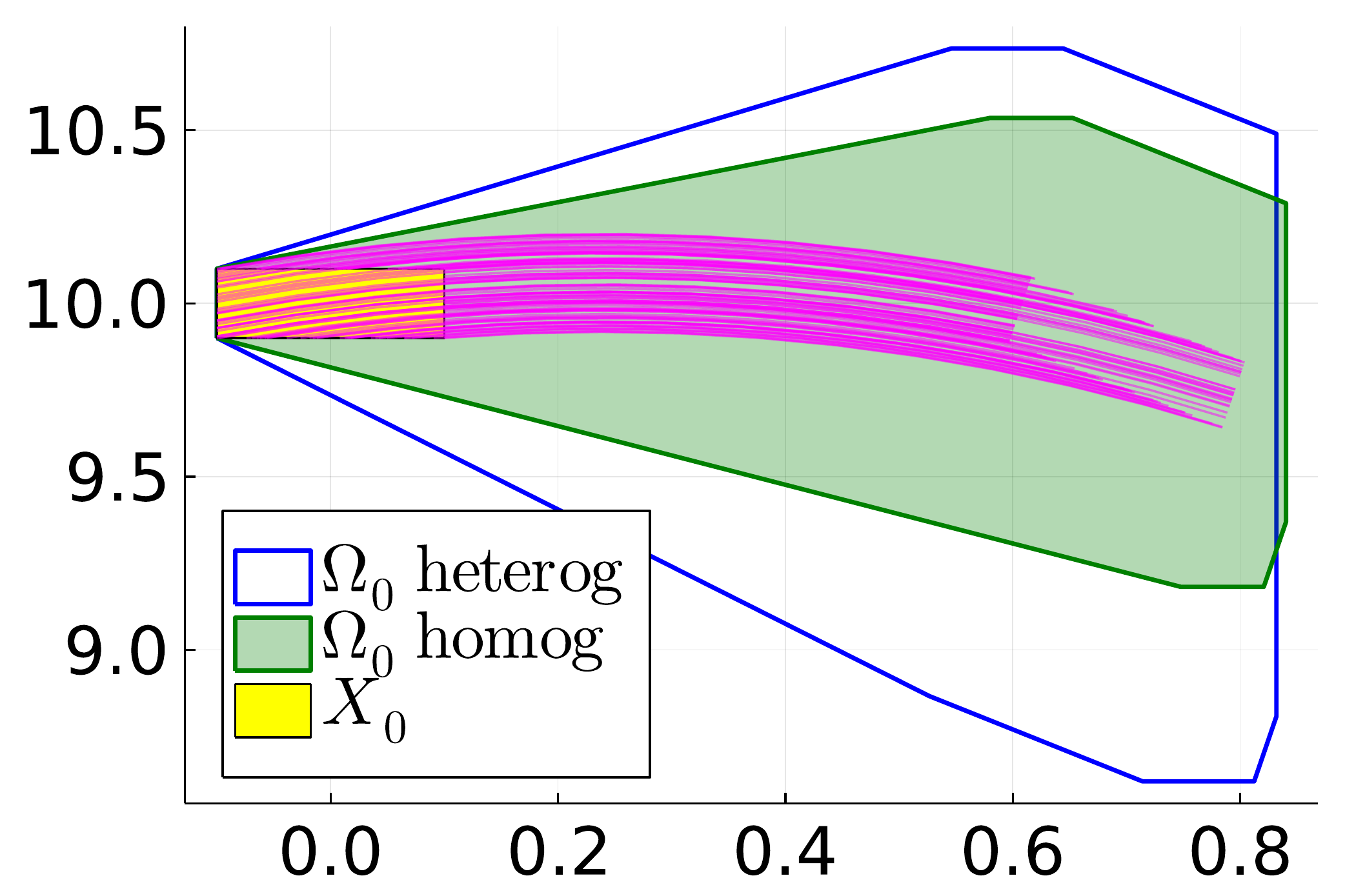}
		\caption{Homogenizing a deterministic system.}
		\label{fig:homogenize_deterministic}
	\end{subfigure}
	\begin{subfigure}[b]{0.49\textwidth}
		\centering
		\includegraphics[width=\linewidth,keepaspectratio]{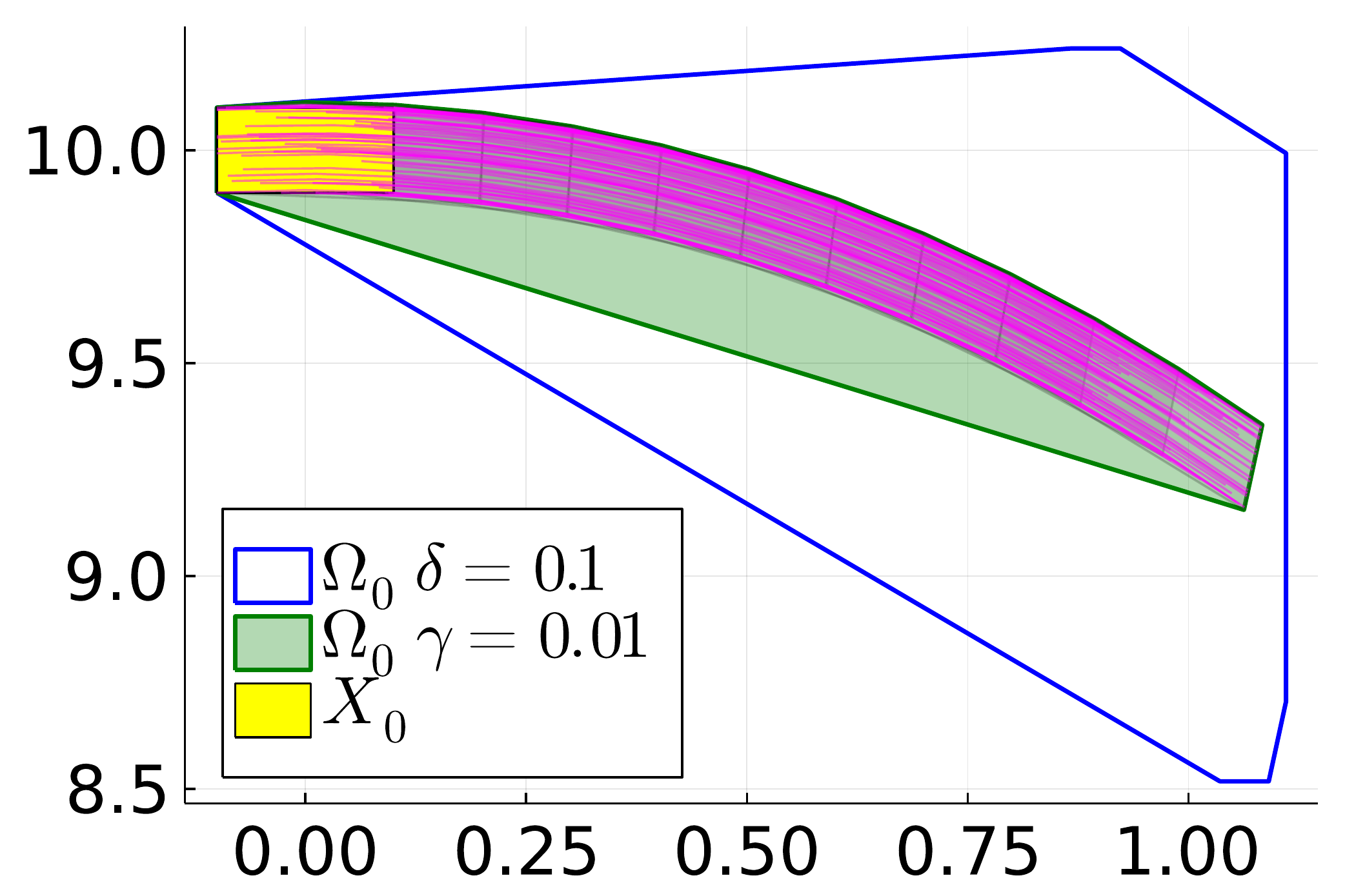}
		\caption{Reduced time step.}
		\label{fig:shrink_delta}
	\end{subfigure}
	\caption{\textbf{Left:} Several trajectories and the sets $\Omega_0$ for a deterministic heterogeneous system (blue) and for the (projected) homogenized system (green).
	\textbf{Right:} Several trajectories and the sets $\Omega_0$ obtained with a time step $\delta = 0.1$ (blue) and obtained by first computing the sets $\Omega_0', \dots, \Omega_9'$ with a time step $\gamma = 0.01$ (light gray, below the trajectories) and then computing the convex hull of their union (green).
	We use the ``forward-only'' method in all cases.}
	\label{fig:homogenize_deterministic_shrink_delta}
\end{figure}

It is easy to see that the set of trajectories of the latter system, after projecting away the auxiliary dimension, is equivalent to the set of trajectories of the original system.
Yet, discretization methods typically yield more precise results for the latter system.
We illustrate this claim in Fig.~\ref{fig:homogenize_deterministic} for the deterministic harmonic oscillator.
The effect for nondeterministic systems is almost identical and is shown in Fig.~\ref{fig:homogenize_nondeterministic} in Appendix~\ref{sec:additional_plots}.

\subsection{Shrinking the time step}

The time step $\delta$ has a large impact on the precision of $\Omega_0$.
Recall that $\Omega_0$ is convex while $\reach_{[0, \delta]}$ is generally not, so we have $\Omega_0 \supseteq \CH(\reach_{[0, \delta]}) \supseteq \reach_{[0, \delta]}$.
The gaps between these three sets grow with larger $\delta$.
Indeed, most methods have the property that they converge to $\reach_{[0, \delta]}$ for $\delta \to 0$.
The reason for not choosing a very small $\delta$ in practice is that the reachability algorithm requires $\lceil \frac{T}{\delta} \rceil$ steps to cover the time horizon $T$.
Hence one needs to find a balance for $\delta$: small enough to obtain a precise $\Omega_0$ and large enough to be efficient later.

Observe that the gap between $\CH(\reach_{[0, \delta]})$ and $\reach_{[0, \delta]}$ only depends on $\delta$.
Say that we fix $\delta$ to keep $\lceil \frac{T}{\delta} \rceil$ in a feasible range.
The best we can do is minimize the gap between $\Omega_0$ and $\CH(\reach_{[0, \delta]})$.
We can achieve this by choosing a positive integer $k$ and a smaller time step $\gamma = \delta / k$, computing the corresponding discretization $\Omega_0'$ and propagating it (using some reachability algorithm) until time horizon $\delta$, which yields sets $\Omega_0', \dots, \Omega_{k-1}'$, and finally constructing $\Omega_0$ from their convex hull $\CH(\Omega_0' \cup \dots \cup \Omega_{k-1}')$.
We illustrate this idea in Fig.~\ref{fig:shrink_delta}.

\section{Efficient implementation}\label{sec:implementation}

As mentioned in Section~\ref{sec:structure}, the different methods can be phrased in the structure of Eq.~\eqref{eq:Omega0_general}.
Exceptions to this rule are the ``first-order zonotope'' method and the ``forward-backward'' method. The former applies an intermediate simplification to the $CH$, while the latter involves a $CH$ in a continuous variable.
In this section we describe how to operate with Eq.~\eqref{eq:Omega0_general} numerically in an efficient way, even in high dimensions.
We show example code for the set library LazySets.jl\footnote{See \href{https://github.com/JuliaReach/LazySets.jl}{github.com/JuliaReach/LazySets.jl}} \cite{LazySets}.
If we allow for a post-processing operation, all the discretization methods can be cast into the same symbolic-numeric framework. We briefly describe the user interface for the implementation of the discretization methods available in the library ReachabilityAnalysis.jl\footnote{See \href{https://github.com/JuliaReach/ReachabilityAnalysis.jl}{github.com/JuliaReach/ReachabilityAnalysis.jl}} in Appendix~\ref{sec:user_interface}.

\subsection{The concept of a lazy set}

The key to an efficient implementation is lazy evaluation, that is, to \emph{delay} computational effort until a result is needed.
For example, we show how the ``first-order d/dt'' method from Section~\ref{ref:firstorder_ddt} can be implemented using LazySets, with the harmonic oscillator from Section~\ref{sec:transformation} as running example.
In this case, the sets $\mathcal{H}$ and $\mathcal{J}$ in Eq.~\eqref{eq:Omega0_general} are the empty set and $\B_\varepsilon$, respectively. The following command defines the set $\Omega_0$ as a lazy representation of Eq.~\eqref{eq:method_homogeneous}.

\begin{minipage}{.93\linewidth}
	\begin{lstlisting}
# first-order d/dt method
julia> Ω₀ = CH(X₀, Φ*X₀) ⊕ B
\end{lstlisting}
\end{minipage}

\noindent%
\begin{minipage}{0.77\linewidth}

In other words, the lazy set operations (Minkowski sum, convex hull, and linear map) are not evaluated.
The execution is instantaneous, while obtaining a concrete representation such as a polyhedron scales with the dimension of the sets.
The operations are internally represented in the form of a tree as shown in the diagram on the right. Further operations such as the support function, conversion, and approximation can be efficiently applied to that symbolic representation \cite{LazySets}.
\end{minipage}
\hfill
\begin{minipage}{0.2\linewidth}
	\vspace*{-3mm}
	\begin{tikzpicture}[level/.style={sibling distance=8mm,level distance=8mm}]
		\node {$\oplus$}
		child {
			node {\CH}
			child {
				node {$\X_0$}
			}
			child {
				node {\texttt{*}}
				child {
					node {$\Phi$}
				}
				child {
					node {$\X_0$}
				}
			}
		}
		child {
			node {$\B$}
		};
	\end{tikzpicture}
\end{minipage}

\subsection{Computation of matrix functions}\label{sec:matrix_functions}

Some discretization methods require special matrix functions such as $\Phi_2$ defined in Eq.~\eqref{eq:Phi2}. If $A$ is not invertible, it can be obtained as the sub-matrix of the exponential of a higher-order matrix \cite{FrehseGDCRLRGDM11}.
However, for large systems (typically $n > 2000$ depending on the sparsity pattern of $A$), such an approach can be prohibitively expensive. Instead, it is possible to use Krylov-subspace methods as discussed in Section~\ref{sec:krylov}, provided that we reformulate the problem as the action of a matrix function over a direction.
To illustrate this point, consider the ``forward-only'' from Section \ref{sec:forward_only}.
Assume that the system is homogeneous and $\X_0$ is a hyperrectangle with center and radius vectors $c, r \in \R^n$ respectively.
A priori, the Krylov method does not apply to Eq.~\eqref{eq:Eplus} because $\Phi_2$ is acting on the \emph{set} $H_\mathit{in} := \boxdot(A^2 \X_0)$.
However, we observe that $H_\mathit{in}$ is a hyperrectangle centered in the origin with radius $r_\mathit{in} = \vert A^2 c\vert + \vert A^2\vert r$. Therefore, it suffices to compute $\vert \Phi_2(\vert A \vert, \delta) \vert r_\mathit{in} = \Phi_2(\vert A \vert, \delta) r_\mathit{in}$ using Krylov methods.

As an application, we consider the discretized heat partial differential equation model used in \cite{ARCH-COMP}, whose details can be found in Appendix~\ref{sec:models}. A run-time comparison is presented in Table~\ref{tab:heat_runtime}. The non-Krylov implementation cannot handle the larger instances because the computation of the matrix exponential runs out memory. The Krylov method additionally offers a significant speedup.

\begin{table}[tb]
	\caption{Average run times (in seconds) for different heat model instances.}
	\label{tab:heat_runtime}
	\centering
	\begin{tabular}{c @{\hspace*{3mm}} | @{\hspace*{3mm}} c @{\hspace*{3mm}} c @{\hspace*{3mm}} c @{\hspace*{3mm}} c}
		\toprule
		Model instance & HEAT01 & HEAT02 & HEAT03 & HEAT04 \\
		\# Mesh points & $5^3$ & $10^3$ & $20^3$ & $50^3$ \\
		\midrule
		Forward &
		\num{0.134} &
		\num{23.6}\phantom{000} &
		-  &
		- \\
		Forward \& Krylov &
		\num{0.001} &
		\num{0.004} &
		\num{0.07} &
		\num{1.367} \\
		\bottomrule
	\end{tabular}
\end{table}

\subsection{Simplification of the set representation}

It is possible to post-process the set $\Omega_0$ with another set that makes it easier to operate with in reachability algorithms. For example, an approximation with an axis-aligned box can be used. The main advantage of such a representation is that the support function can be computed efficiently.
This is shown in the following comparison when computing the support function of $\Omega_0$ resp.\ its box approximation along direction $d = (1, 1)^T$.
The box approximation of $\Omega_0$ has approximately the same support value, but the computation is 13x faster.

\begin{minipage}{.93\linewidth}
	\begin{lstlisting}
julia> boxΩ₀ = box_approximation(Ω₀);

julia> d = [1.0, 1.0];

julia> ρ(d, Ω₀)
10.328776223585699

julia> ρ(d, boxΩ₀)
10.35390370135613

julia> @btime ρ($d, $Ω₀)
97.762 ns (1 allocation: 80 bytes)
  
julia> @btime ρ($d, $boxΩ₀)
7.209 ns (0 allocations: 0 bytes)
\end{lstlisting}
\end{minipage}

\section{Experimental evaluation}\label{sec:evaluation}

In this section we evaluate the discretization methods from Section~\ref{sec:methods} in two experiments.\footnote{The scripts are available at~\cite{RE}.}
In the first experiment, we visually compare the sets $\Omega_0$ for variants of the harmonic oscillator (see Section~\ref{sec:transformation}).
In the second experiment, we evaluate the methods on representative models while varying the time step $\delta$.
Next we describe the experimental setup, show the results, and finally discuss them.

\subsection{Setup}

We implemented the different methods (called here d/dt, Zonotope, Correction hull, First-order, Forward/backward, and Forward) in JuliaReach \cite{BogomolovFFPS19}.
For matrix-exponential functionality we use \cite{ExponentialUtilities}.
For the correction-hull method we use the truncation order $p = 4$ (higher orders led to numerical errors for the biggest model).
To plot results for the forward-backward method, one needs to choose a set of directions to evaluate the support function.
In the plots we choose $30$ uniform directions from the unit circle; we choose this high precision to show the theoretical power of the method, even if in practice this is rarely required.
We obtained the results on a notebook with a 2.20\,GHz CPU and 8\,GB RAM.

\subsection{Visual evaluation of varying parameters}\label{sec:visual_evaluation}

\begin{figure}[p]
	\centering
	\includegraphics[width=.49\linewidth,keepaspectratio]{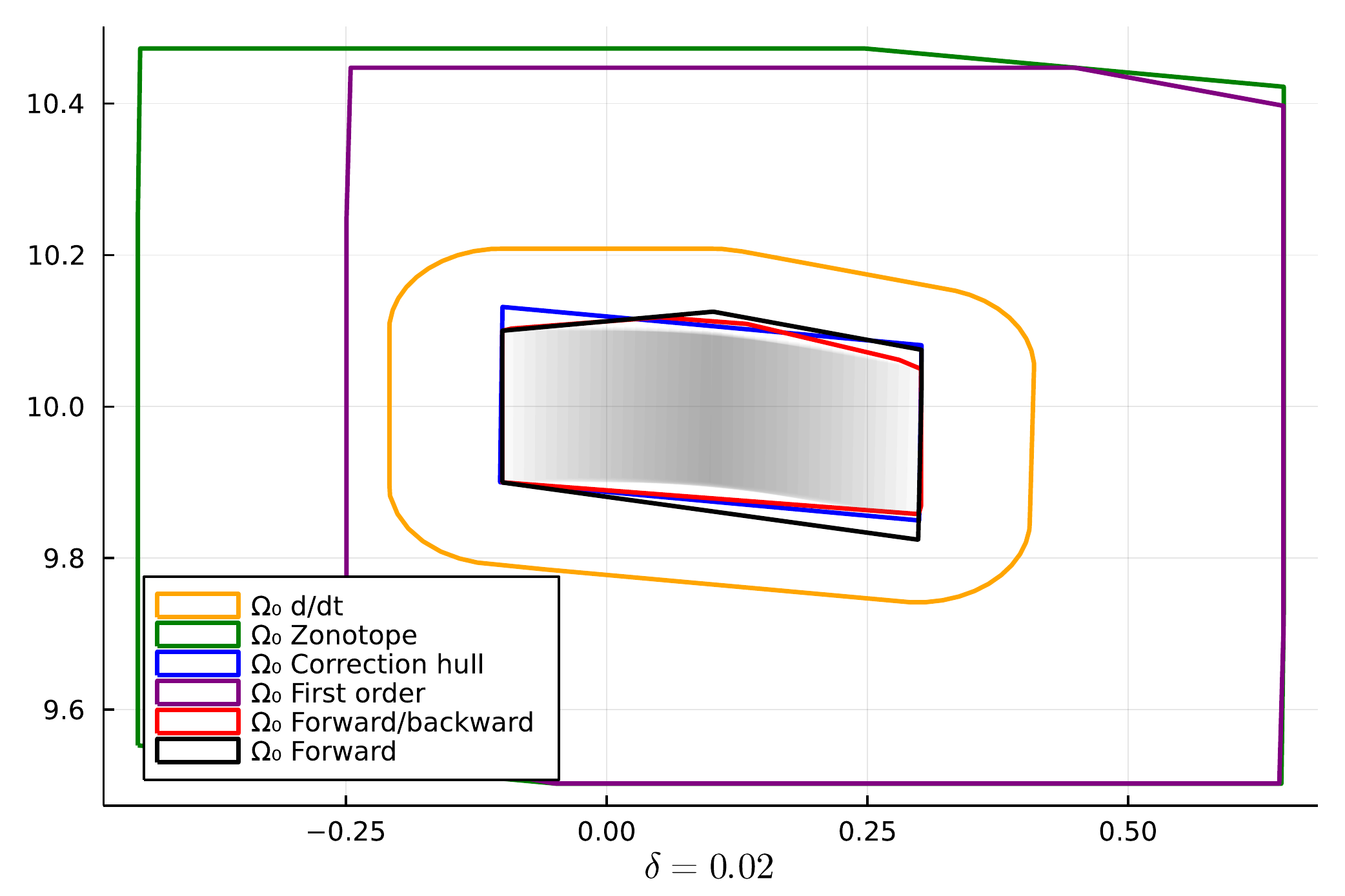}
	\includegraphics[width=.49\linewidth,keepaspectratio]{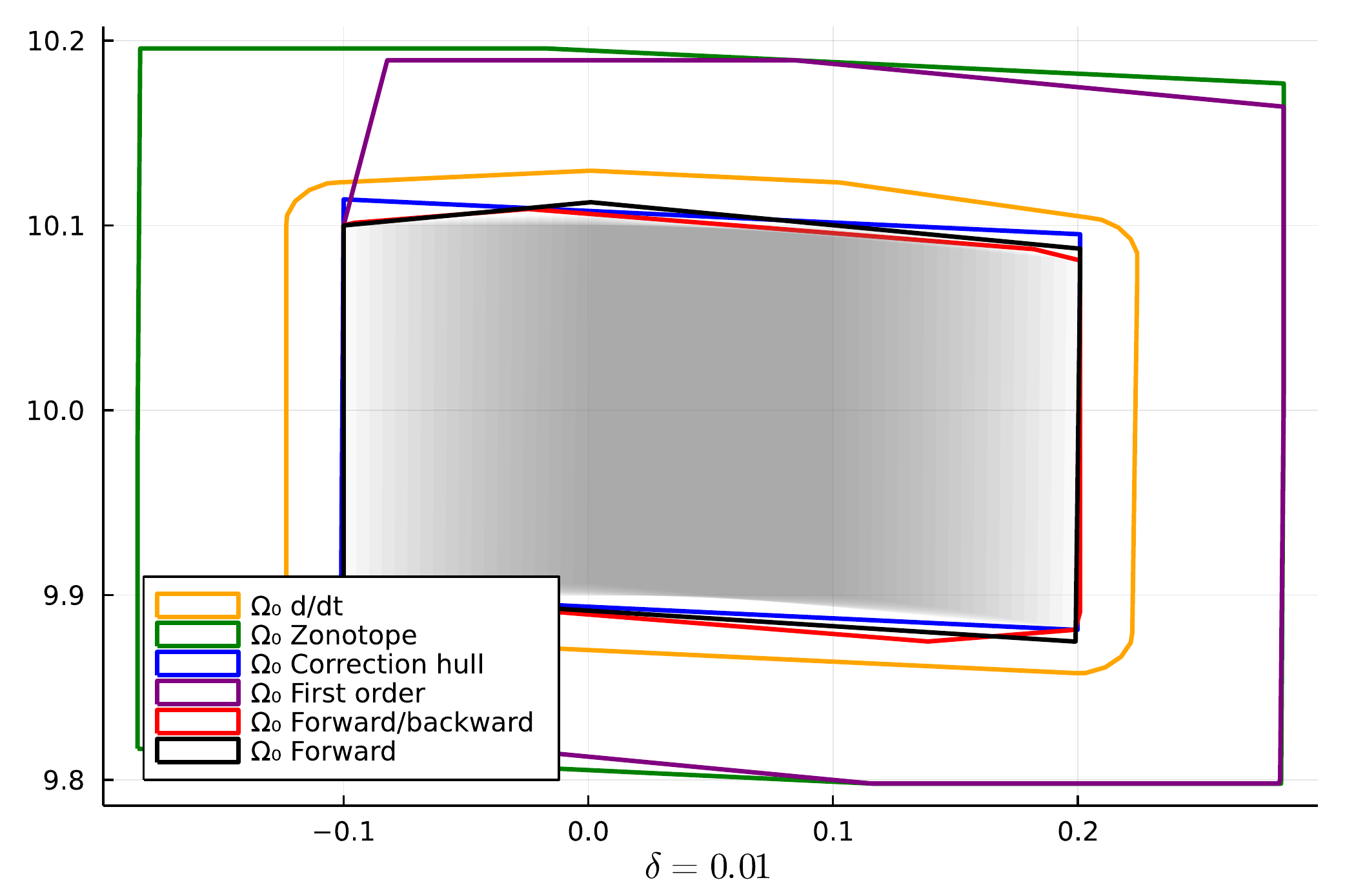}
	\includegraphics[width=.49\linewidth,keepaspectratio]{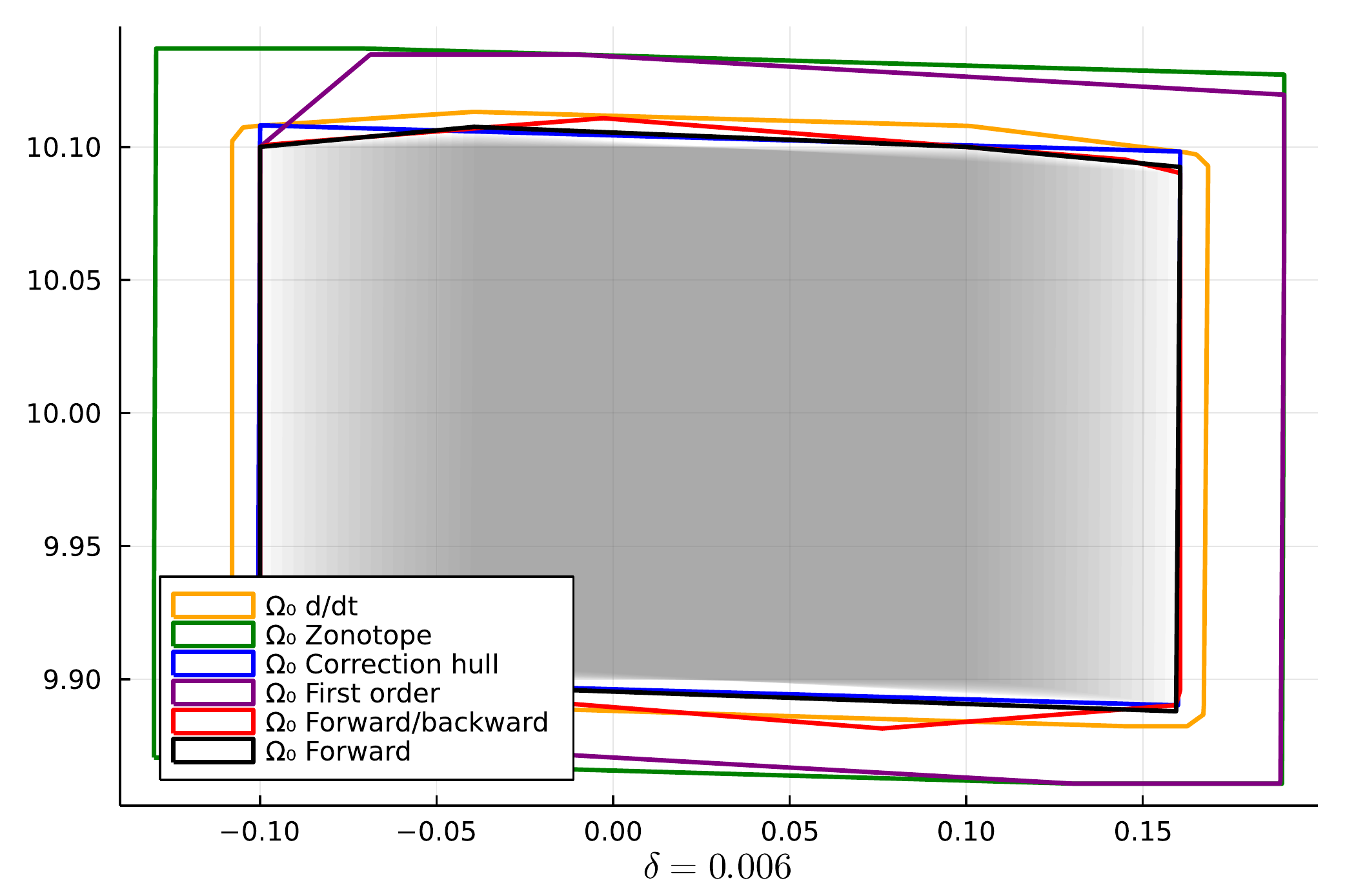}
	\includegraphics[width=.49\linewidth,keepaspectratio]{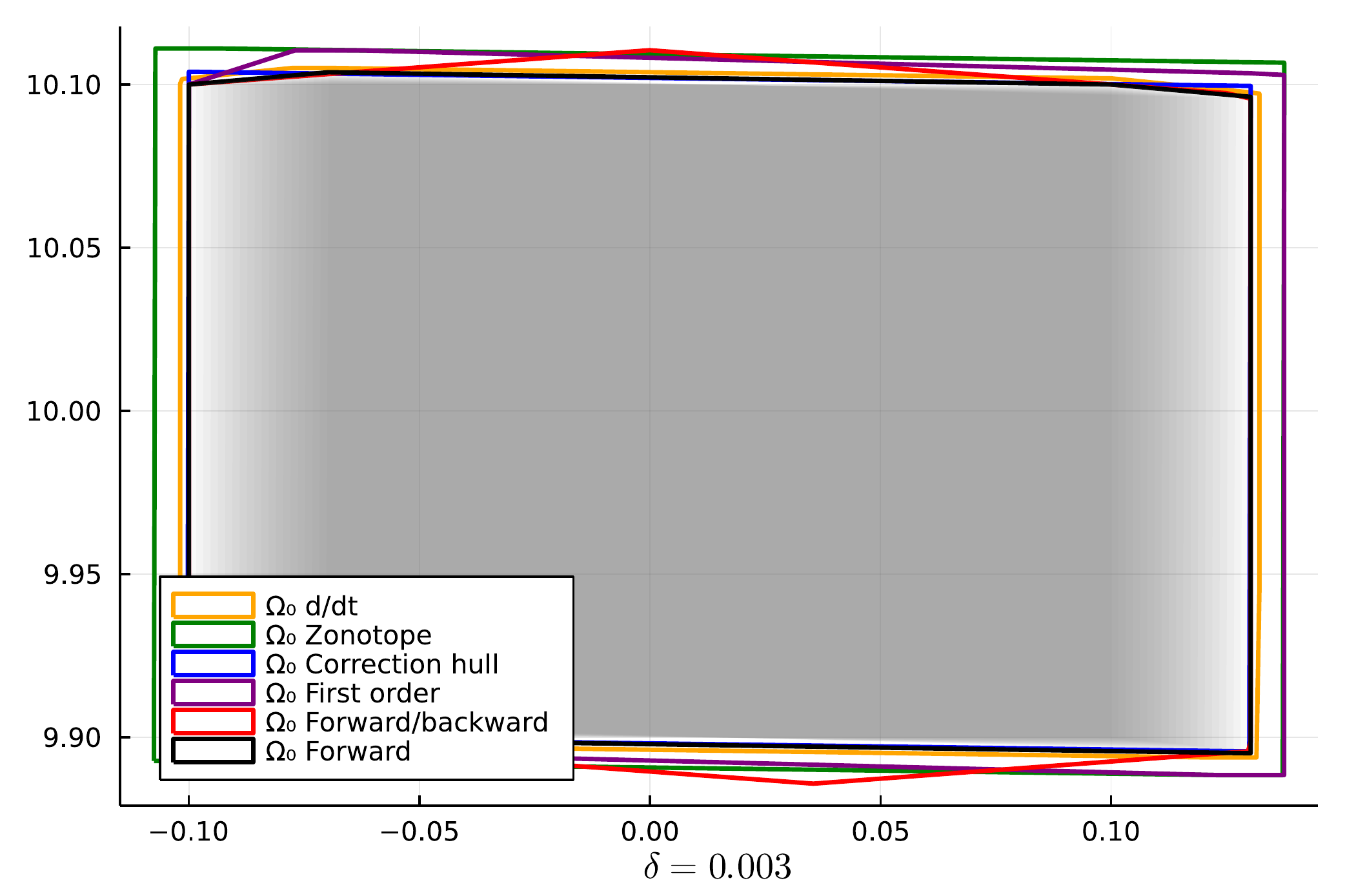}
	\caption{The sets $\Omega_0$ obtained with different methods for varying values of $\delta$.
		In gray we show $\reach_t$ for uniform $t \in [0, \delta]$.}
	\label{fig:experiment_delta}
\end{figure}

\begin{figure}[p]
	\centering
	\includegraphics[width=.49\linewidth,keepaspectratio]{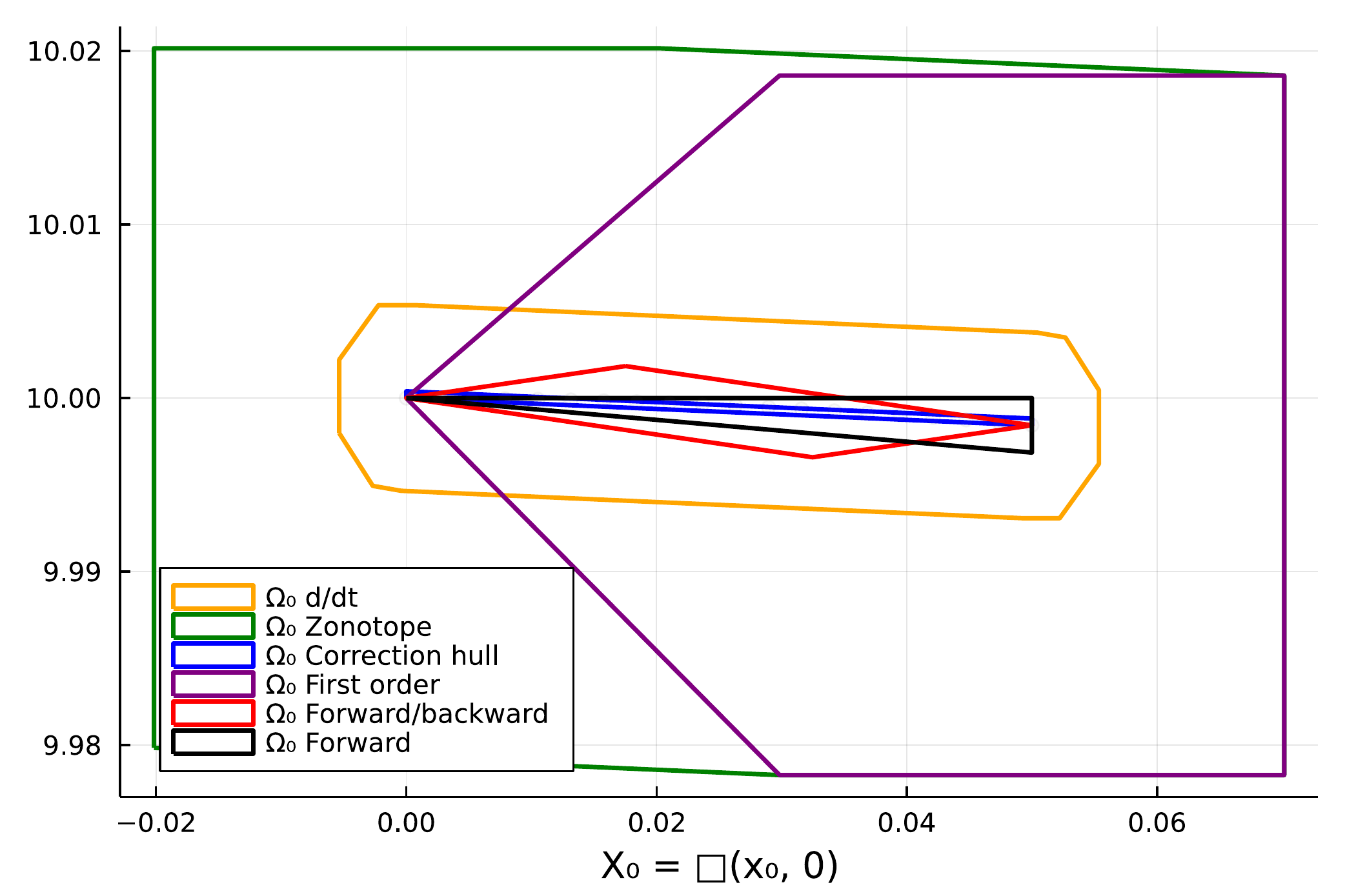}
	\includegraphics[width=.49\linewidth,keepaspectratio]{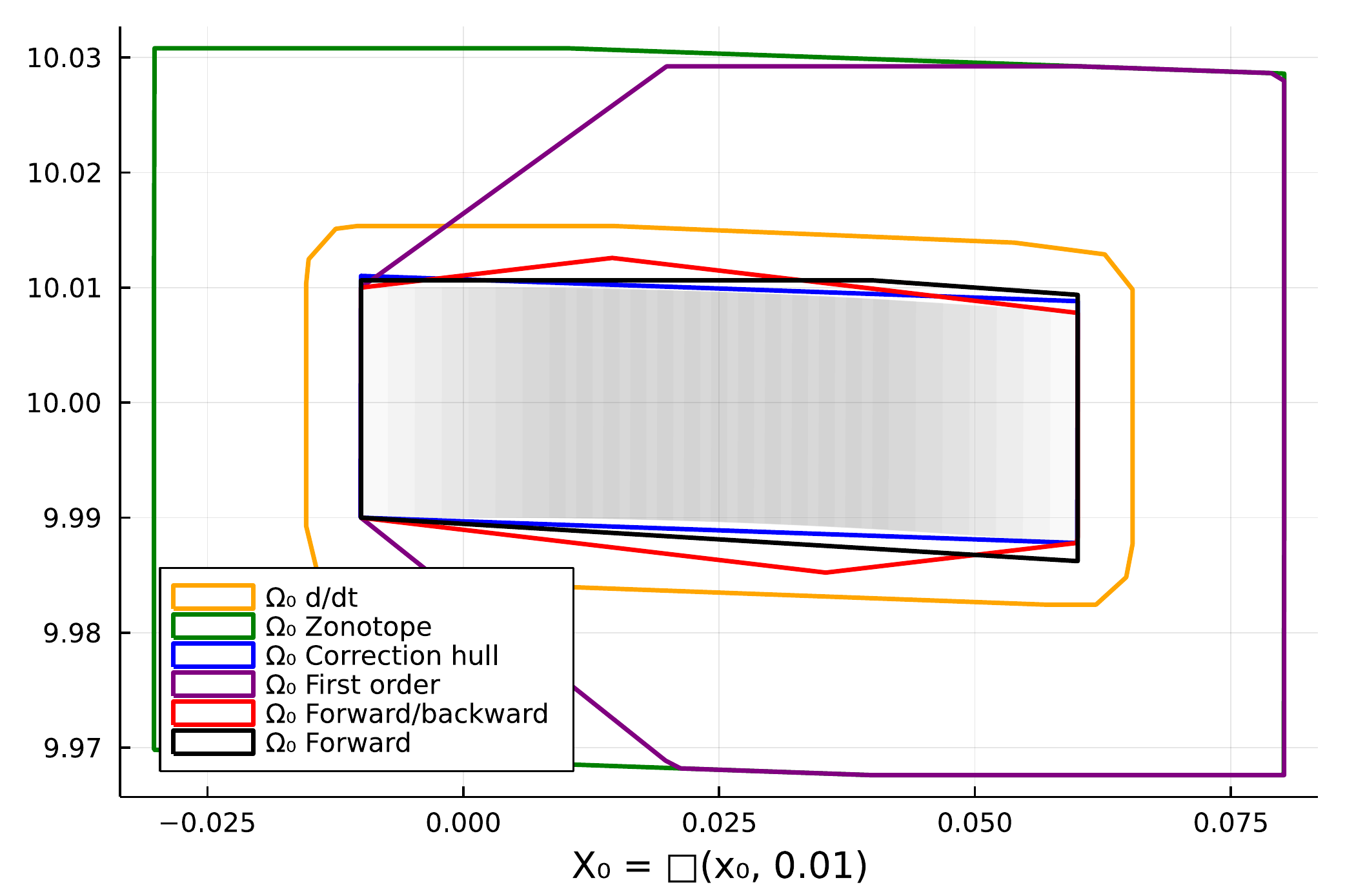}
	\includegraphics[width=.49\linewidth,keepaspectratio]{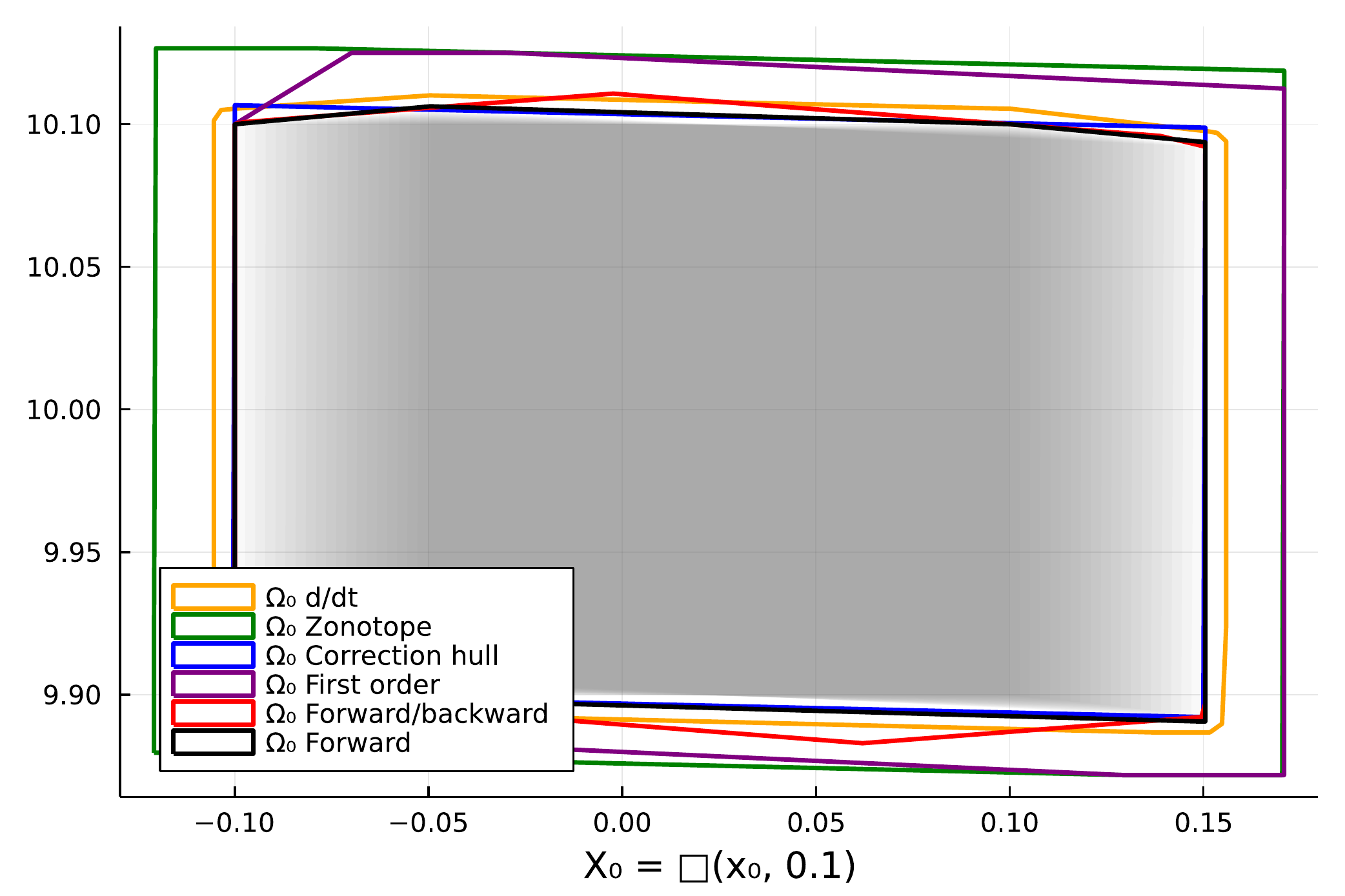}
	\includegraphics[width=.49\linewidth,keepaspectratio]{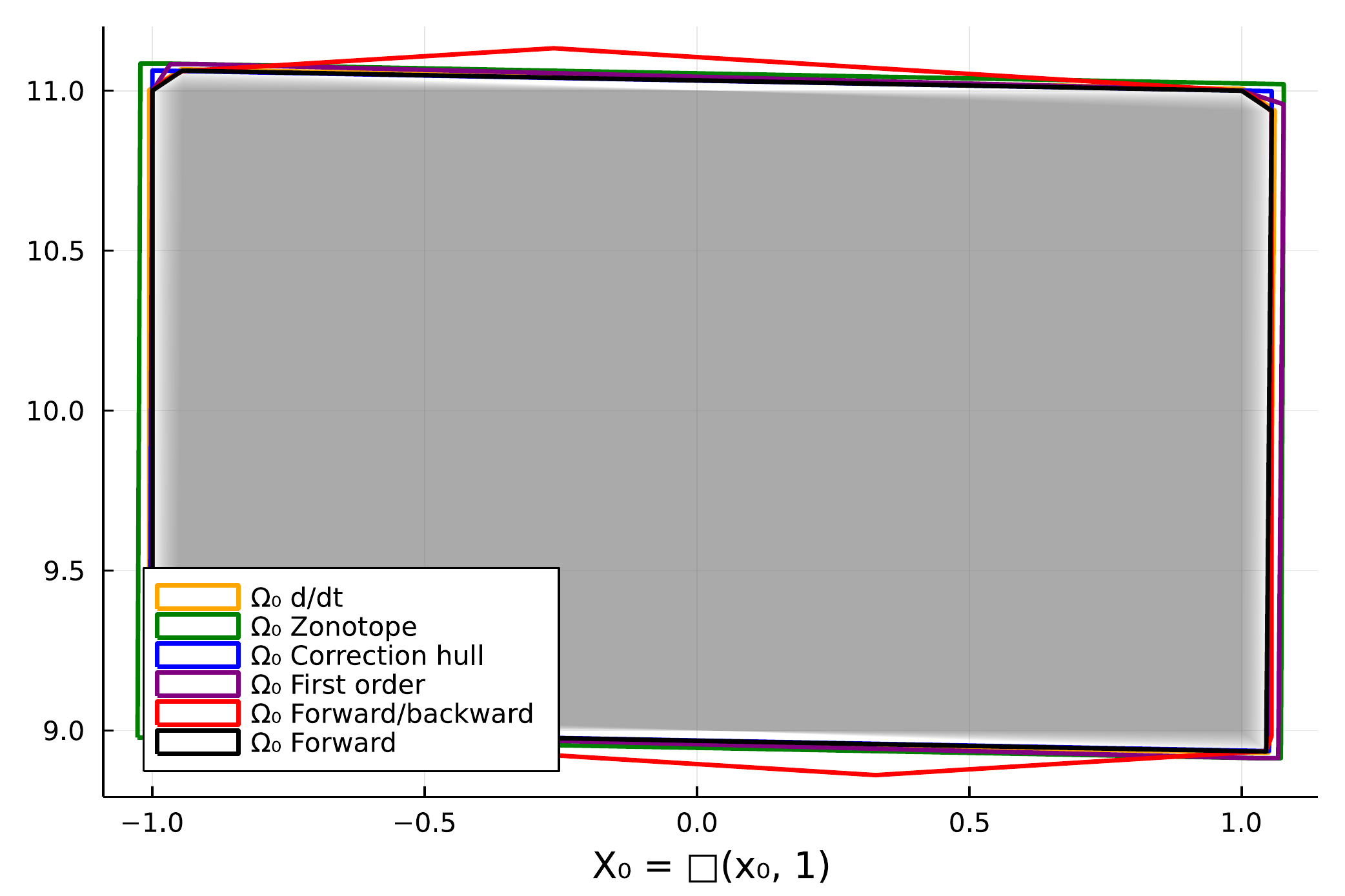}
	\caption{The sets $\Omega_0$ obtained with different methods for $\delta = 0.005$ and varying sets $\X_0$.
	In gray we show $\reach_t$ for uniform $t \in [0, \delta]$.}
	\label{fig:experiment_X0}
\end{figure}

\begin{figure}[t]
	\centering
	\includegraphics[width=.49\linewidth,keepaspectratio]{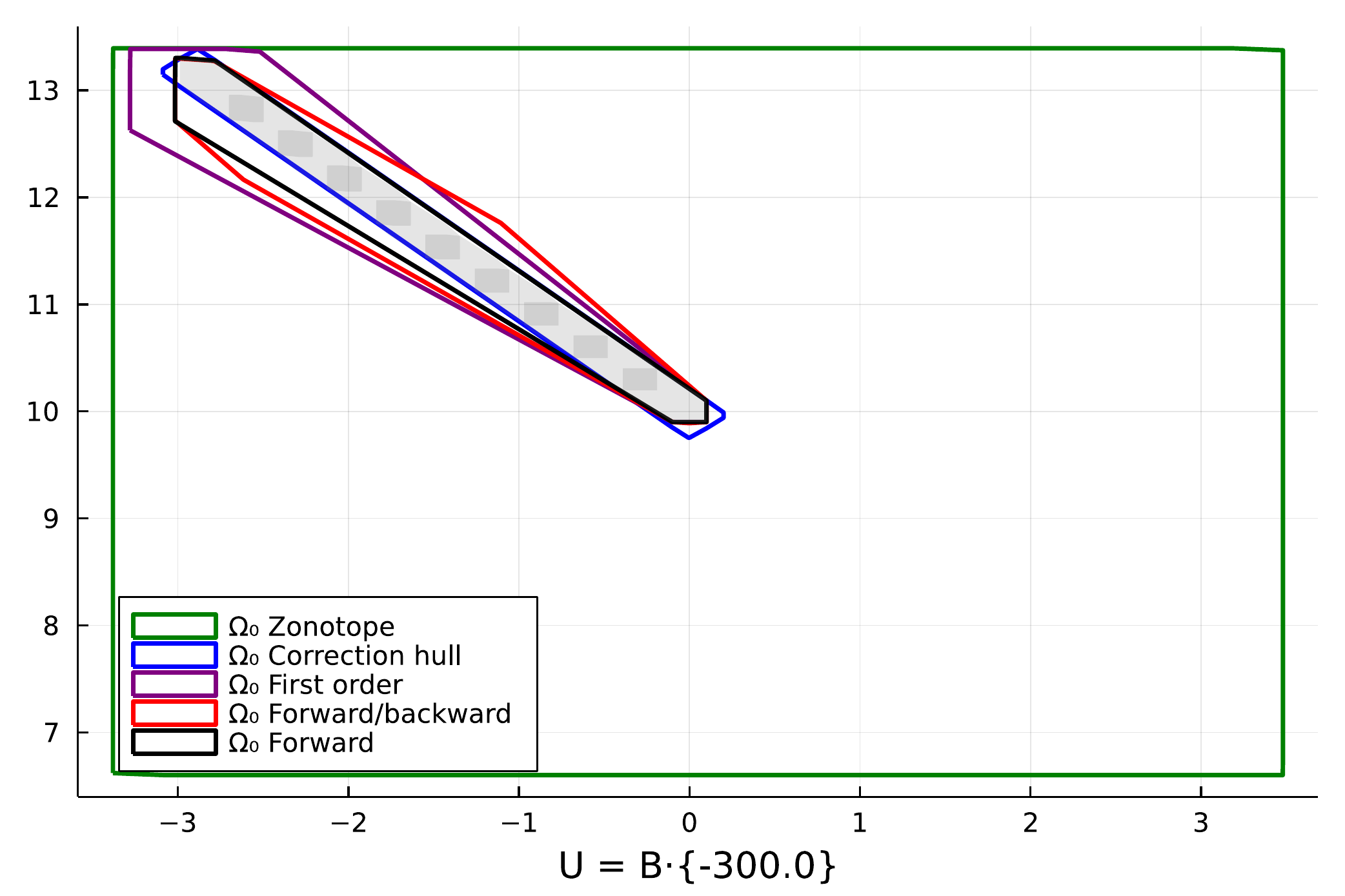}
	\includegraphics[width=.49\linewidth,keepaspectratio]{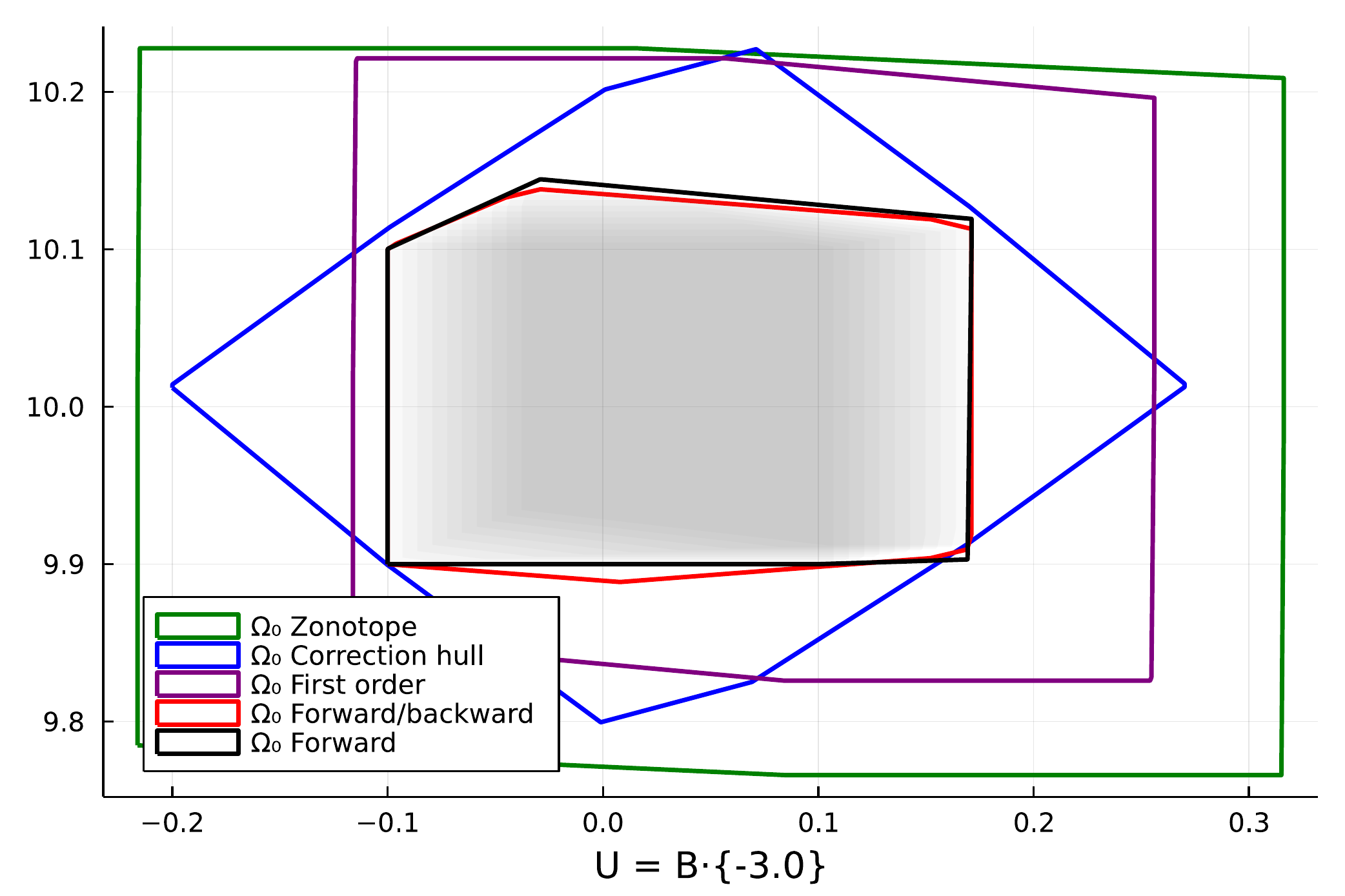}
	\includegraphics[width=.49\linewidth,keepaspectratio]{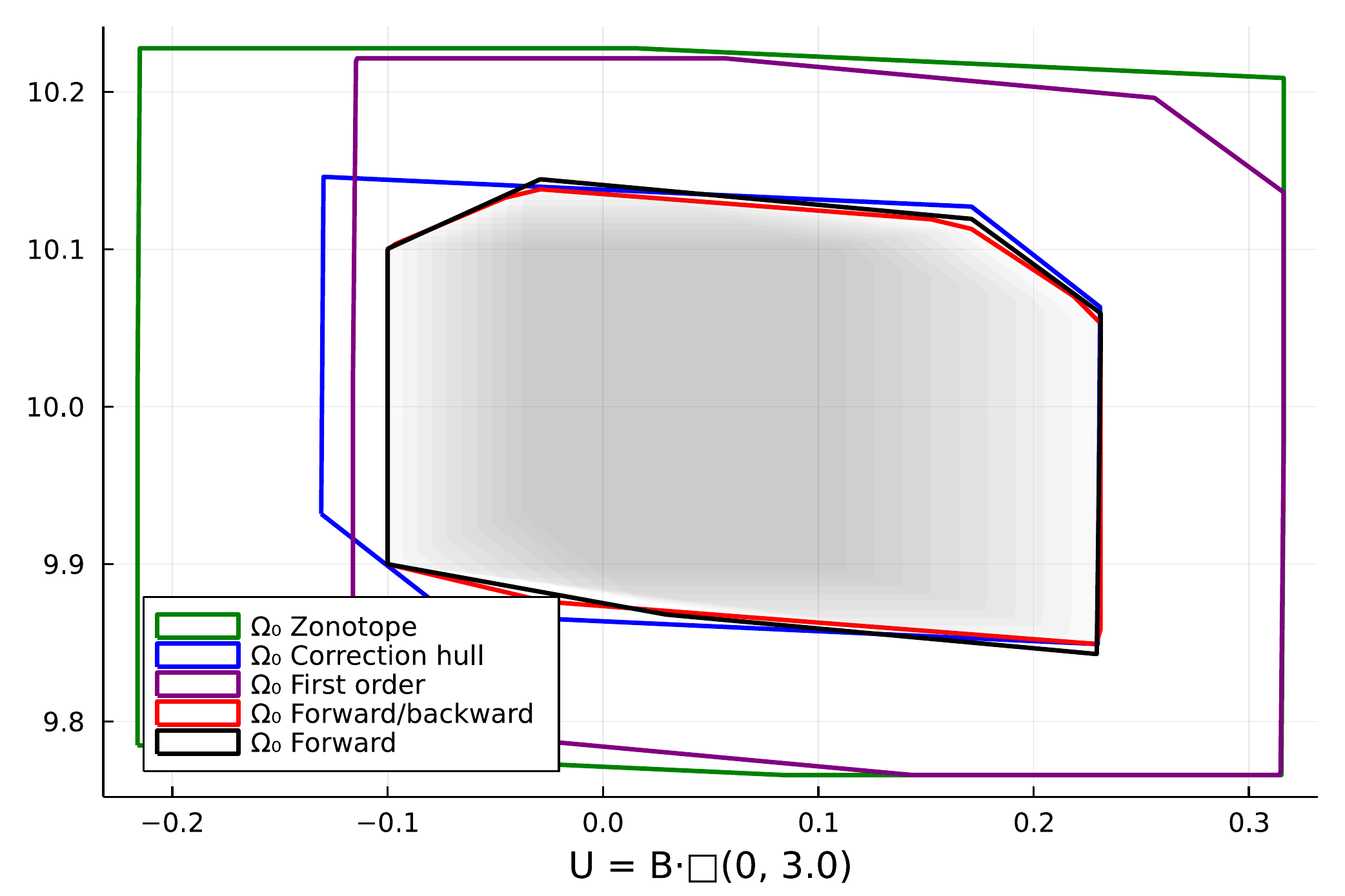}
	\includegraphics[width=.49\linewidth,keepaspectratio]{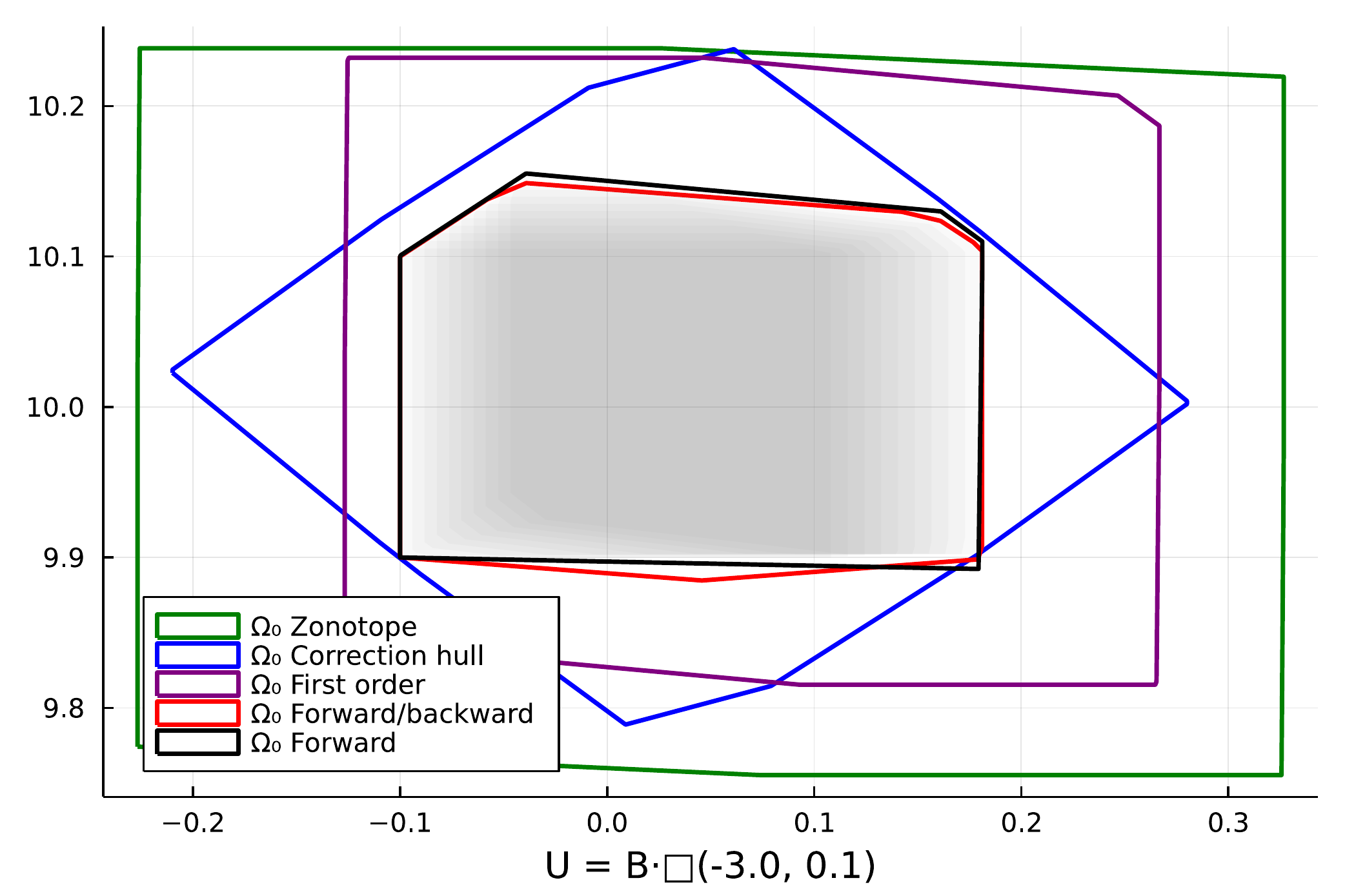}
	\caption{The sets $\Omega_0$ obtained with different methods for $\delta = 0.01$ and varying sets $\U$.
	In gray we show $\reach_t$ for uniform $t \in [0, \delta]$.}
	\label{fig:experiment_U}
\end{figure}

We evaluate the methods on the harmonic oscillator for three different analysis and model parameters.
In Fig.~\ref{fig:experiment_delta} we vary the step size $\delta$.
In Fig.~\ref{fig:experiment_X0} we vary the size of the initial set $\X_0$.
In Fig.~\ref{fig:experiment_U} we vary the size of the input domain \U.
The plots also show a tight approximation of the true reachable states.
For homogeneous systems, the analytic solution at time $t$ can be computed ($e^{At} \X_0$) and we show several sets $\reach_t$ (for uniformly chosen time points $t$ from $[0, \delta]$) instead.

\subsection{Quantitative evaluation by scaling $\delta$}\label{sec:benchmark}

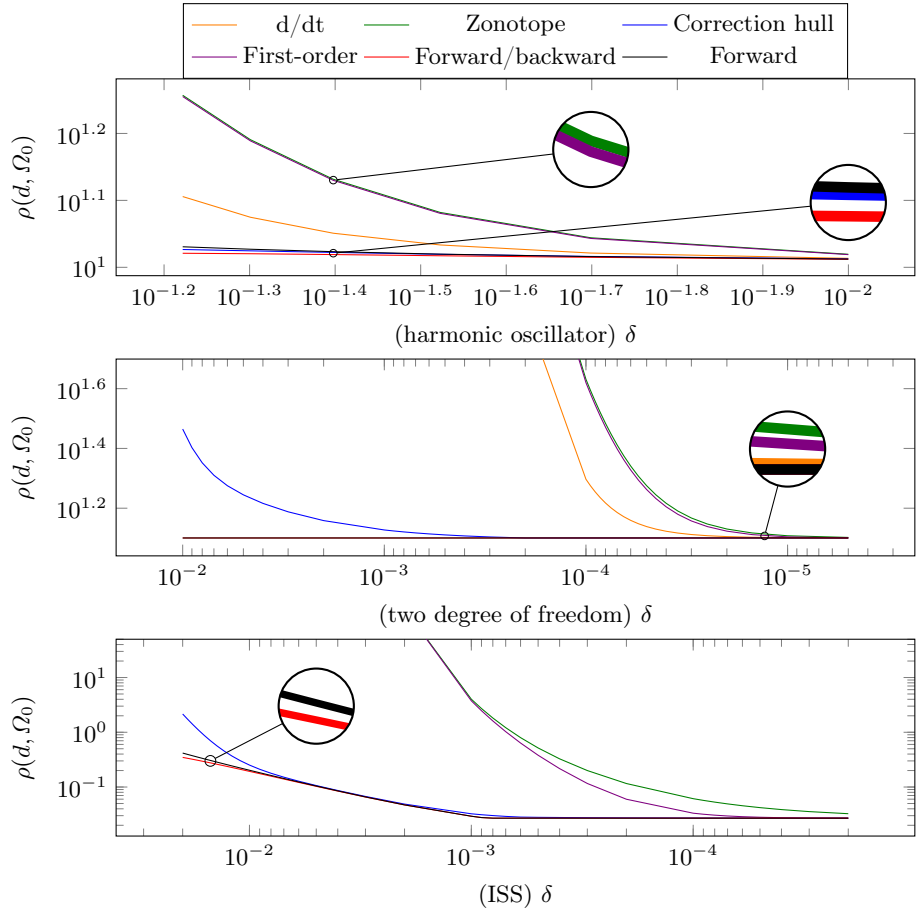
\begin{figure}[t!]
	\centering
	\begin{tikzpicture}[spy using outlines={circle, magnification=10, connect spies}]
		\begin{loglogaxis}[height=42mm,width=\linewidth,no markers,xlabel={(harmonic oscillator) $\delta$},ylabel={$\rho(d, \Omega_0)$},x dir=reverse,legend columns=3,legend style={at={(0.5,1)},anchor=south}]
			\addplot +[COLOR1] table[x index=0, y index=1] {delta_support_oscillator.dat};
			\addplot +[COLOR2] table[x index=0, y index=2] {delta_support_oscillator.dat};
			\addplot +[COLOR3] table[x index=0, y index=3] {delta_support_oscillator.dat};
			\addplot +[COLOR4] table[x index=0, y index=4] {delta_support_oscillator.dat};
			\addplot +[COLOR5] table[x index=0, y index=5] {delta_support_oscillator.dat};
			\addplot +[COLOR6,solid] table[x index=0, y index=6] {delta_support_oscillator.dat};
			\coordinate (spypoint1) at (axis cs:0.04,10.5);
			\coordinate (magnifyglass1) at (axis cs:0.01,12.5);
			\coordinate (spypoint2) at (axis cs:0.04,13.5);
			\coordinate (magnifyglass2) at (axis cs:0.02,15);
			\addlegendentry{d/dt}
			\addlegendentry{Zonotope}
			\addlegendentry{Correction hull}
			\addlegendentry{First-order}
			\addlegendentry{Forward/backward}
			\addlegendentry{Forward}
		\end{loglogaxis}
		\spy [size=10mm] on (spypoint1) in node[fill=white] at (magnifyglass1);
		\spy [size=10mm] on (spypoint2) in node[fill=white] at (magnifyglass2);
	\end{tikzpicture}
	\begin{tikzpicture}[spy using outlines={circle, magnification=10, connect spies}]
		\begin{loglogaxis}[height=42mm,width=\linewidth,ymax=50,no markers,xlabel={(two degree of freedom) $\delta$},ylabel={$\rho(d, \Omega_0)$},x dir=reverse
		]
			\addplot +[COLOR1] table[x index=0, y index=1] {delta_support_freedom.dat};
			\addplot +[COLOR2] table[x index=0, y index=2] {delta_support_freedom.dat};
			\addplot +[COLOR3] table[x index=0, y index=3] {delta_support_freedom.dat};
			\addplot +[COLOR4] table[x index=0, y index=4] {delta_support_freedom.dat};
			\addplot +[COLOR5] table[x index=0, y index=5] {delta_support_freedom.dat};
			\addplot +[COLOR6,solid] table[x index=0, y index=6] {delta_support_freedom.dat};
			\coordinate (spypoint1) at (axis cs:0.000013,12.8);
			\coordinate (magnifyglass1) at (axis cs:1e-5,25);
		\end{loglogaxis}
		\spy [size=10mm] on (spypoint1) in node[fill=white] at (magnifyglass1);
	\end{tikzpicture}
	\begin{tikzpicture}[spy using outlines={circle, magnification=7, connect spies}]
		\begin{loglogaxis}[height=42mm,width=\linewidth,ymax=50,no markers,xlabel={(ISS) $\delta$},ylabel={$\rho(d, \Omega_0)$},x dir=reverse
		]
			\addplot +[COLOR2] table[x index=0, y index=1] {delta_support_iss.dat};
			\addplot +[COLOR3] table[x index=0, y index=2] {delta_support_iss.dat};
			\addplot +[COLOR4] table[x index=0, y index=3] {delta_support_iss.dat};
			\addplot +[COLOR5] table[x index=0, y index=4] {delta_support_iss.dat};
			\addplot +[COLOR6,solid] table[x index=0, y index=5] {delta_support_iss.dat};
			\coordinate (spypoint) at (axis cs:1.5e-2,0.3);
			\coordinate (magnifyglass) at (axis cs:0.005,3);
		\end{loglogaxis}
		\spy [size=10mm] on (spypoint) in node[fill=white] at (magnifyglass);
	\end{tikzpicture}
	\vspace*{-7mm}
	\caption{Benchmark results with the graphs of $\rho(d, \Omega_0)$ (log axes).
	In the second plot, the methods ``forward-backward'' and ``forward-only'' yield identical results.
	In the third plot, the ``d/dt'' method is not applicable.
	}
	\label{fig:benchmark_delta}
\end{figure}

\pgfkeys{/pgf/fpu,/pgf/fpu/output format=fixed}  

\begin{filecontents*}{delta_time_oscillator.dat}
delta	d-dt	Zonotope	Correction-hull	First-order	Forward-backward	Forward
0.01	0.006793243243243242	0.012716378378378376	0.11420037837837838	0.0037847297297297303	6.50190827027027	0.02630237837837838
0.02	0.015238054054054052	0.026751837837837838	0.25870894594594585	0.008289702702702705	6.395956675675675	0.030968918918918925
0.03	0.015247783783783784	0.026348567567567566	0.25911294594594597	0.008265189189189191	6.364605702702702	0.03023764864864865
0.04	0.015342783783783785	0.025918837837837834	0.2568906216216216	0.008333702702702702	6.77689991891892	0.029923378378378378
0.05	0.013238459459459461	0.023479864864864867	0.22987556756756758	0.006606243243243243	6.46639035135135	0.030215135135135143
0.06	0.015918594594594593	0.02507227027027027	0.2639014594594595	0.008215324324324324	6.872763135135136	0.03583781081081081
\end{filecontents*}
\csvreader[before reading=\def\total{0},separator=tab]{delta_time_oscillator.dat}{delta=\delta}{%
\pgfmathsetmacro{\total}{\total+1}%
}
\csvreader[before reading=\def\Osumddt{0},separator=tab]{delta_time_oscillator.dat}{d-dt=\ddt}{%
\pgfmathsetmacro{\Osumddt}{\Osumddt+\ddt}%
}
\pgfmathsetmacro{\Osumddt}{\Osumddt/\total}
\csvreader[before reading=\def\Osumzono{0},separator=tab]{delta_time_oscillator.dat}{Zonotope=\zono}{%
\pgfmathsetmacro{\Osumzono}{\Osumzono+\zono}%
}
\pgfmathsetmacro{\Osumzono}{\Osumzono/\total}
\csvreader[before reading=\def\Osumcorrhull{0},separator=tab]{delta_time_oscillator.dat}{Correction-hull=\corrhull}{%
\pgfmathsetmacro{\Osumcorrhull}{\Osumcorrhull+\corrhull}%
}
\pgfmathsetmacro{\Osumcorrhull}{\Osumcorrhull/\total}
\csvreader[before reading=\def\Osumfirstorder{0},separator=tab]{delta_time_oscillator.dat}{First-order=\firstorder}{%
\pgfmathsetmacro{\Osumfirstorder}{\Osumfirstorder+\firstorder}%
}
\pgfmathsetmacro{\Osumfirstorder}{\Osumfirstorder/\total}
\csvreader[before reading=\def\Osumforwback{0},separator=tab]{delta_time_oscillator.dat}{Forward-backward=\forwback}{%
\pgfmathsetmacro{\Osumforwback}{\Osumforwback+\forwback}
}
\pgfmathsetmacro{\Osumforwback}{\Osumforwback/\total}
\csvreader[before reading=\def\Osumforw{0},separator=tab]{delta_time_oscillator.dat}{Forward=\forw}{%
\pgfmathsetmacro{\Osumforw}{\Osumforw+\forw}%
}
\pgfmathsetmacro{\Osumforw}{\Osumforw/\total}

\pgfkeys{/pgf/fpu=false}

\pgfkeys{/pgf/fpu,/pgf/fpu/output format=fixed}  

\begin{filecontents*}{delta_time_freedom.dat}
delta	d-dt	Zonotope	Correction-hull	First-order	Forward-backward	Forward
5.0e-6	0.026339405405405407	0.03205543243243244	0.37579440540540543	0.01306864864864865	5.869009864864865	0.054335540540540544
1.0e-5	0.02667886486486487	0.039092756756756754	0.4945607297297297	0.014660054054054053	5.933539702702703	0.05399897297297297
1.5e-5	0.026587081081081085	0.038614	0.48568518918918924	0.01294072972972973	5.7543111351351355	0.05402529729729731
2.0e-5	0.02713859459459459	0.03960813513513512	0.513563054054054	0.013269918918918917	5.897210837837836	0.055492243243243256
2.5e-5	0.02578308108108108	0.03837191891891892	0.4507225135135135	0.013673378378378377	5.775761864864865	0.05787297297297298
3.0e-5	0.025735837837837835	0.03712351351351352	0.34321932432432434	0.012703216216216215	6.3473793783783785	0.049325405405405406
3.5e-5	0.026364918918918918	0.03158948648648648	0.454703054054054	0.013509135135135136	5.954942027027028	0.05274424324324324
4.0e-5	0.0249367027027027	0.039174243243243236	0.45374427027027026	0.012661756756756758	5.700654216216216	0.05246602702702703
4.5e-5	0.025815945945945943	0.03578027027027027	0.49994783783783786	0.012763270270270272	5.738022945945946	0.05920881081081081
5.0e-5	0.028252594594594595	0.03499656756756757	0.4411214864864865	0.01367945945945946	6.078552702702702	0.053086162162162165
5.5e-5	0.02588816216216216	0.038795243243243245	0.4402313243243243	0.013937648648648649	5.79623445945946	0.11554305405405406
6.0e-5	0.05867227027027027	0.08567797297297298	0.9311637027027028	0.03047286486486486	6.680209432432433	0.1605291891891892
6.5e-5	0.08032613513513513	0.10694197297297299	0.7289837837837836	0.015119621621621624	6.203135513513514	0.055580297297297285
7.0e-5	0.03603291891891892	0.041224351351351364	0.5071238108108108	0.01392754054054054	6.0683305945945945	0.06017997297297298
7.5e-5	0.02736481081081081	0.03855513513513514	0.5091961621621621	0.014728945945945947	5.502870702702703	0.062440432432432434
8.0e-5	0.02874981081081081	0.041511162162162156	0.5457367027027028	0.014756270270270273	5.278544135135134	0.053231648648648655
8.5e-5	0.030143	0.03691751351351351	0.43124135135135133	0.012894243243243245	5.169598918918918	0.05420002702702703
9.0e-5	0.02699164864864865	0.03983786486486486	0.503214783783784	0.012643783783783785	5.072768459459458	0.05247875675675676
9.5e-5	0.027000540540540543	0.03571627027027027	0.4305	0.013080027027027027	5.663321297297298	0.054858270270270265
0.0001	0.028935486486486486	0.036261675675675675	0.5334846486486485	0.013342432432432435	6.012948702702704	0.051953567567567575
0.0002	0.060504081081081074	0.3103571081081082	1.269702135135135	0.014886378378378378	5.8589983243243235	0.06492145945945946
0.0003	0.029297621621621618	0.041325891891891904	0.504732162162162	0.01703975675675676	6.890274027027027	0.057058081081081076
0.0004	0.02698859459459459	0.040341108108108106	0.45937708108108105	0.013864162162162161	6.48295108108108	0.0644248108108108
0.0005	0.026354459459459453	0.03786864864864866	0.539553891891892	0.01370145945945946	6.481423135135135	0.06035929729729729
0.0006	0.15085613513513513	0.03969016216216217	0.4656176216216216	0.013930756756756753	6.910121891891893	0.053380945945945946
0.0007	0.028228189189189196	0.03513543243243243	0.447294054054054	0.012684000000000003	6.625782378378379	0.053861729729729725
0.0008	0.025825054054054053	0.0343212972972973	0.39155205405405413	0.014772540540540544	6.527693999999999	0.05748818918918918
0.0009	0.026617513513513512	0.0368341081081081	0.5158506486486486	0.013929513513513512	6.507744675675675	0.056840702702702715
0.001	0.026631648648648643	0.03669375675675675	0.473368972972973	0.013872081081081084	6.386431027027028	0.06428054054054054
0.002	0.032982594594594586	0.03811921621621621	0.5233875945945945	0.014993486486486485	6.544646378378379	0.06132148648648647
0.003	0.027509432432432434	0.038370918918918924	0.4510236756756756	0.015442945945945946	6.849018756756757	0.06234956756756758
0.004	0.027811189189189185	0.03702489189189189	0.39139302702702705	0.016557432432432434	7.241521918918918	0.062146729729729726
0.005	0.03355259459459458	0.043406162162162164	0.5377448918918919	0.015406621621621625	7.19170743243243	0.06234072972972972
0.006	0.03179548648648649	0.04013648648648648	0.4704975945945946	0.015275864864864864	6.253163297297297	0.06254127027027027
0.007	0.029688864864864863	0.04170481081081082	0.4389089189189189	0.015070675675675677	6.224318810810811	0.06433227027027028
0.008	0.0287227027027027	0.04023732432432432	0.4614031081081081	0.018366540540540544	6.256615648648649	0.06623975675675677
0.009	0.029146891891891895	0.040101189189189194	0.5269142972972972	0.016447081081081078	6.292678243243243	0.06788797297297296
0.01	0.02916883783783784	0.04548843243243243	0.4893754054054053	0.016708864864864868	6.249152540540541	0.0675087027027027
\end{filecontents*}
\csvreader[before reading=\def\total{0},separator=tab]{delta_time_freedom.dat}{delta=\delta}{%
\pgfmathsetmacro{\total}{\total+1}%
}
\csvreader[before reading=\def\Fsumddt{0},separator=tab]{delta_time_freedom.dat}{d-dt=\ddt}{%
\pgfmathsetmacro{\Fsumddt}{\Fsumddt+\ddt}%
}
\pgfmathsetmacro{\Fsumddt}{\Fsumddt/\total}
\csvreader[before reading=\def\Fsumzono{0},separator=tab]{delta_time_freedom.dat}{Zonotope=\zono}{%
\pgfmathsetmacro{\Fsumzono}{\Fsumzono+\zono}%
}
\pgfmathsetmacro{\Fsumzono}{\Fsumzono/\total}
\csvreader[before reading=\def\Fsumcorrhull{0},separator=tab]{delta_time_freedom.dat}{Correction-hull=\corrhull}{%
\pgfmathsetmacro{\Fsumcorrhull}{\Fsumcorrhull+\corrhull}%
}
\pgfmathsetmacro{\Fsumcorrhull}{\Fsumcorrhull/\total}
\csvreader[before reading=\def\Fsumfirstorder{0},separator=tab]{delta_time_freedom.dat}{First-order=\firstorder}{%
\pgfmathsetmacro{\Fsumfirstorder}{\Fsumfirstorder+\firstorder}%
}
\pgfmathsetmacro{\Fsumfirstorder}{\Fsumfirstorder/\total}
\csvreader[before reading=\def\Fsumforwback{0},separator=tab]{delta_time_freedom.dat}{Forward-backward=\forwback}{%
\pgfmathsetmacro{\Fsumforwback}{\Fsumforwback+\forwback}
}
\pgfmathsetmacro{\Fsumforwback}{\Fsumforwback/\total}
\csvreader[before reading=\def\Fsumforw{0},separator=tab]{delta_time_freedom.dat}{Forward=\forw}{%
\pgfmathsetmacro{\Fsumforw}{\Fsumforw+\forw}%
}
\pgfmathsetmacro{\Fsumforw}{\Fsumforw/\total}

\pgfkeys{/pgf/fpu=false}

\pgfkeys{/pgf/fpu,/pgf/fpu/output format=fixed}  

\begin{filecontents*}{delta_time_iss.dat}
delta	Zonotope	Correction-hull	First-order	Forward-backward	Forward
2.0e-5	12.337652	4109.381582999999	12.307309	587.391424	617.8194410000001
3.0e-5	91.296455	5175.88255	11.495285	566.8504149999999	371.381631
4.0e-5	18.588381000000002	4068.2274740000003	7.781540000000001	533.006836	363.55642700000004
5.0e-5	12.212966999999999	4030.081726	8.000555	535.543711	365.697887
6.0e-5	13.284996	4237.3729969999995	9.715459000000001	567.400538	452.717394
7.0e-5	20.182213	5815.4696460000005	117.97122300000001	766.4749210000001	370.40275199999996
8.0e-5	13.498237000000001	4993.0981	9.411581	581.5462699999999	477.90456500000005
9.0e-5	23.876402	5308.167876	115.68211300000002	1096.566127	775.2618749999999
0.0001	168.028006	4846.052745999999	10.52682	595.720009	383.98269400000004
0.0002	18.039212	6474.510943	81.82816799999999	1174.1251809999999	688.475806
0.0003	20.216327000000003	6389.686189999999	9.551069	637.458628	415.55942799999997
0.0004	31.553606000000002	4358.206099000001	9.224089	554.9433919999999	392.958112
0.0005	21.128929	4235.397088000001	12.550971000000002	583.308442	386.334241
0.0006	18.090469	4343.620713	13.115439999999998	565.516413	436.693474
0.0007	18.019650000000002	4425.970235	9.820720999999999	574.183332	413.933904
0.0008	18.240597	4376.71121	9.265661999999999	579.665794	401.66486499999996
0.0009	23.295662999999998	4317.775259000001	9.981992	576.1397	471.62533099999996
0.001	27.458457	5390.452959	12.603824000000001	579.8382490000001	396.318067
0.002	119.538298	4551.894071	9.333681	564.131574	402.490924
0.003	18.480856	4610.190540999999	20.42503	605.41276	398.064611
0.004	19.424057	4178.468119	12.22497	582.284317	416.50825900000007
0.005	18.720553	4122.266295	11.417642	675.302851	429.524706
0.006	18.894771000000002	4434.58903	18.660453	630.6341889999999	869.3384920000001
0.007	92.182396	5189.603395000001	12.009547	636.615295	469.18917899999997
0.008	21.560032999999997	5230.292705999999	14.191215	645.672681	427.899855
0.009	19.430746	5183.865412	12.935592	649.3926540000001	740.3091039999999
0.01	61.092592	6412.091448000001	213.05492	1352.076723	531.3952350000001
0.011	18.494529	4291.591457	11.578726	593.638661	414.676355
0.012	18.113438000000002	4275.257344000001	22.241802	612.631302	415.728233
0.013	18.261129	4241.947529	10.618293000000001	666.274551	591.4675830000001
0.014	19.879247	4081.020921	10.935262000000002	584.5270059999999	437.932742
0.015	45.316436	4119.355204	11.945537000000002	639.8686130000001	475.895885
0.016	33.463026000000006	4270.253387	11.789704	646.753033	453.128817
0.017	21.395177	4183.825634999999	11.815904	600.53066	448.956953
0.018	23.440060000000003	4236.979119	12.527115	622.882395	460.068066
0.019	20.239666	4702.917186	30.664565	730.195208	456.779656
0.02	43.328953	4738.243371	30.324354	644.854248	626.005247
\end{filecontents*}
\csvreader[before reading=\def\total{0},separator=tab]{delta_time_iss.dat}{delta=\delta}{%
\pgfmathsetmacro{\total}{\total+1}%
}
\csvreader[before reading=\def\Isumzono{0},separator=tab]{delta_time_iss.dat}{Zonotope=\zono}{%
\pgfmathsetmacro{\Isumzono}{\Isumzono+\zono}
}
\pgfmathsetmacro{\Isumzono}{\Isumzono/\total}
\csvreader[before reading=\def\Isumcorrhull{0},separator=tab]{delta_time_iss.dat}{Correction-hull=\corrhull}{%
\pgfmathsetmacro{\Isumcorrhull}{\Isumcorrhull+\corrhull}%
}
\pgfmathsetmacro{\Isumcorrhull}{\Isumcorrhull/\total}
\csvreader[before reading=\def\Isumfirstorder{0},separator=tab]{delta_time_iss.dat}{First-order=\firstorder}{%
\pgfmathsetmacro{\Isumfirstorder}{\Isumfirstorder+\firstorder}%
}
\pgfmathsetmacro{\Isumfirstorder}{\Isumfirstorder/\total}
\csvreader[before reading=\def\Isumforwback{0},separator=tab]{delta_time_iss.dat}{Forward-backward=\forwback}{%
\pgfmathsetmacro{\Isumforwback}{\Isumforwback+\forwback}
}
\pgfmathsetmacro{\Isumforwback}{\Isumforwback/\total}
\csvreader[before reading=\def\Isumforw{0},separator=tab]{delta_time_iss.dat}{Forward=\forw}{%
\pgfmathsetmacro{\Isumforw}{\Isumforw+\forw}%
}
\pgfmathsetmacro{\Isumforw}{\Isumforw/\total}

\pgfkeys{/pgf/fpu=false}

\begin{table}[tb]
	\caption{Average run times (in milliseconds) for the different methods.}
	\label{tab:runtime}
	\centering
	\begin{tabular}{c | r @{\hspace*{1.8mm}} r @{\hspace*{1.8mm}} r @{\hspace*{1.8mm}} r @{\hspace*{1.8mm}} r @{\hspace*{1.8mm}} r}
		\toprule
		Model & \mc{d/dt} & \mc{Zonotope} & \mc{Correction hull} & \mc{First-order} & \mc{Fwd/bwd} & \mc{Forward} \\
		\midrule
		Oscillator &
		\pgfmathroundtozerofill{\Osumddt} \num{\pgfmathresult} &
		\pgfmathroundtozerofill{\Osumzono} \num{\pgfmathresult} &
		\pgfmathroundtozerofill{\Osumcorrhull} \num{\pgfmathresult} &
		\pgfmathroundtozerofill{\Osumfirstorder} \num{\pgfmathresult} &
		\pgfmathroundtozerofill{\Osumforwback} \num{\pgfmathresult} &
		\pgfmathroundtozerofill{\Osumforw} \num{\pgfmathresult} \\
		TDoF &
		\pgfmathroundtozerofill{\Fsumddt} \num{\pgfmathresult} &
		\pgfmathroundtozerofill{\Fsumzono} \num{\pgfmathresult} &
		\pgfmathroundtozerofill{\Fsumcorrhull} \num{\pgfmathresult} &
		\pgfmathroundtozerofill{\Fsumfirstorder} \num{\pgfmathresult} &
		\pgfmathroundtozerofill{\Fsumforwback} \num{\pgfmathresult} &
		\pgfmathroundtozerofill{\Fsumforw} \num{\pgfmathresult} \\
		ISS &
		\mc{--} &
		\pgfmathroundtozerofill{\Isumzono} \num{\pgfmathresult} &
		\pgfmathroundtozerofill{\Isumcorrhull} \num{\pgfmathresult} &
		\pgfmathroundtozerofill{\Isumfirstorder} \num{\pgfmathresult} &
		\pgfmathroundtozerofill{\Isumforwback} \num{\pgfmathresult} &
		\pgfmathroundtozerofill{\Isumforw} \num{\pgfmathresult} \\
		\bottomrule
	\end{tabular}
\end{table}

\let\delta\deltatmp  

We run a quantitative analysis on the harmonic oscillator and two other models described in Appendix~\ref{sec:models}.
The latter represent challenging model classes:
The second model, a two-degree-of-freedom system, has a system matrix $A$ of large norm; this shows the corresponding sensitivity of the methods.
The third model, representing a docking maneuver at the International Space Station (ISS), is high-dimensional; this shows the scalability of the methods.
For comparing the precision, we vary the time step $\delta$ and compute the support function $\rho(d, \Omega_0)$ of the sets $\Omega_0$ in direction $d = (1, \dots, 1)^T$.
The results are shown in Fig.~\ref{fig:benchmark_delta}.
Average run times (which are independent of $\delta$) are given in Table~\ref{tab:runtime}.
Note that $\Omega_0$ is computed lazily except for the ``zonotope'' and the ``correction-hull'' methods.

\subsection{Summary}

We generally observe that the first-oder methods (d/dt, Zonotope, First-order) yield coarser results and are more sensitive to the different model characteristics.
In particular, for the two degree of freedom, a very small time step $10^{-5}$ is required to obtain precise results.
This shows the sensitivity to the norm of $A$.
The other three methods yield higher and similar precision, although the ``correction-hull'' method, which it is restricted to zonotopes, is generally incomparable (see for instance the first plot in Fig.~\ref{fig:experiment_U}).
The ``forward-backward'' method is typically the most precise but also the most expensive method; recall that we computed a lazy representation of $\Omega_0$ here and in a reachability application one needs to evaluate the support function of $\Omega_0$ in multiple directions.
The forward-only method is a simplification that yields a good compromise.
The ``correction-hull'' method is much slower than the other methods for the largest model (ISS), but it computes a concrete zonotope.
(Zonotopes enable an efficient ``second stage'' reachability analysis, which we ignore here.)
Since this method is designed for interval matrices, it may be possible to devise a more efficient scalar variant.

\section{Conclusion}

In this article we have studied six methods for conservative time discretization.
We have discussed potential ways to improve their output in practice and how to efficiently implement them.
Our empirical evaluation shows that the methods have different characteristics.
In particular, methods based on a first-order approximation are generally less precise than the other methods.
In the future we plan to perform a similar study of the full-fledged reachability algorithms using the insights gained in this study to use a precise discretization model.

\subsubsection{Acknowledgments}

This research was partly supported by DIREC - Digital Research Centre Denmark and the Villum Investigator Grant S4OS.

\bibliographystyle{splncs04}
\bibliography{bibliography}

\clearpage
\appendix

\section{Model definitions}\label{sec:models}

\subsection{Heat3D}

The linear Heat3D model is obtained from a spatial discretization of a partial differential equation for heat transfer in three dimensions. Originally presented in \cite{HanK06} for two dimensions and later extended to three dimensions, it was used as a benchmark example for reachability analysis \cite{BakTJ19,ARCH-COMP}. The model dimension is scalable and here we consider four different instances of increasing complexity, which are labeled \texttt{HEAT0x}, for grids of size $n \times n \times n$ mesh points, i.e., the associated ODEs are $n^3$-dimensional. The goal is to find the maximum temperature reached at the center of the spatially discretized domain, where one of its edges is initially heated. Since each mesh point corresponds to a given direction, it is sufficient to compute the support function along the center of the mesh. Furthermore, the set of initial states is a hyperrectangle contained in $[0, 1]^3$ and the matrix $A$ is hermitian.

For the experiment described in Section~\ref{sec:matrix_functions} we use Krylov subspace dimension $m = 30$ for instance \texttt{HEAT01} and $m = 100$ for the rest (see \cite{koskela2015approximating} for details).

\subsection{Two degree of freedom}

We consider a two-degree-of-freedom model from \cite[Chapter 9]{hughes2012finite}.
The model has characteristics that are typical of large systems, containing both low-frequency and high-frequency components. It is given by
\begin{equation}\label{eq:2dof}
	M\ddot{y}(t) + Ky(t) = 0,\qquad y(t) = [y_1(t), y_2(t)]^T
\end{equation}
where the mass ($M$) and stiffness ($K$) matrices are respectively
\begin{equation}
M = \begin{pmatrix}
m_1 & 0 \\ 0 & m_2
\end{pmatrix},\qquad K = \begin{pmatrix}
k_1 + k_2 & -k_2 \\ -k_2 & k_2
\end{pmatrix}.
\end{equation}
Eq.~\eqref{eq:2dof} is brought to first-order form of Eq.~\eqref{eq:LTI_general} by introducing the variable $x(t) = [y_1(t), y_2(t), \dot{y}_1(t), \dot{y}_2(t)]^T$, from which we obtain
\begin{equation}\label{eq:2dof_Amatrix}
	\dot{x}(t) = A x(t), \quad \text{where } A = \begin{pmatrix}
		0 & 0 & 1 & 0 \\
		0 & 0 & 0 & 1 \\
		-\frac{k_1+k_2}{m_1} & \frac{k_2}{m_1} & 0 & 0 \\
		\frac{k_2}{m_2} & -\frac{k_2}{m_2} & 0 & 0
	\end{pmatrix}.
\end{equation}
The initial states are the box centered at $[1, 10, 0, 0]$ with radius $[0.1, 0.5, 0.5, 0.5]$.
The numerical values for the parameters are $m_1 = m_2 = k_2 = 1$ and $k_1 = 10{,}000$.

\subsection{ISS}

The ISS (International Space Station) model was originally presented in \cite{ChahlaouiD05} and later proposed as a benchmark example for order-reduction techniques \cite{TranNJ16} and reachability analysis \cite{ARCH-COMP}.
It is a structural model of the component 1R (Russian service module) and models the vibration during the docking maneuver of a spacecraft.
There are $270$ state variables and three nondeterministic inputs.
The matrices $A$ and $e^A$ are sparse (${>}\,99\%$ sparsity) with $\Vert A\Vert_\infty \approx 3763$, $\Vert \X_0 \Vert_\infty = 10^{-4}$, and $\Vert\U\Vert_\infty \approx 0.98$.

\section{User interface}\label{sec:user_interface}

The methods described in Section \ref{sec:methods} have been implemented in the library \href{http://github.com/JuliaReach/ReachabilityAnalysis.jl}{ReachabilityAnalysis.jl} in the Julia programming language.
The implementation offers a user interface that allows to choose different algorithms that can be used through the \texttt{discretize} function.
It receives three arguments: an initial-value problem of the form \eqref{eq:LTI_general} (special cases included), the step size $\delta$, and an algorithm (which defaults to the \texttt{Forward} method).
The function returns a discrete-time system whose set of initial states is $\Omega_0$, and can be obtained from the examples below with the command \texttt{initial\_states(Di)}.

\begin{minipage}{.93\linewidth}
	\begin{lstlisting}
# load reachability library
using ReachabilityAnalysis

# state matrix
A = [0.0 1; -4π 0]

# step-size
δ = 1e-2

# initial states (infinity-norm ball)
X0 = BallInf([0.0, 10.0], 0.1)

# initial-value problem
P = @ivp(x' = A*x, x(0) ∈ X0)

# use default discretization (Forward-only)
D1 = discretize(P, δ)

# use first-order ddt discretization
D2 = discretize(P, δ, FirstOrderddt())

# use correction hull with default order (p = 10)
D3 = discretize(P, δ, CorrectionHull())

# use correction hull with specified order (p = 15)
D4 = discretize(P, δ, CorrectionHull(order=15))
\end{lstlisting}
\end{minipage}

Finally, we mention that the library offers a \texttt{solve} function that first calls \texttt{discretize} and then uses the desired set propagation algorithm to compute the set of states reachable until a given time horizon. These steps correspond to the first (resp. second) stages according the notation in Section~\ref{sec:introduction}. We refer to the online documentation of \href{http://github.com/JuliaReach/ReachabilityAnalysis.jl}{ReachabilityAnalysis.jl} for further examples.

\section{Underapproximation}\label{sec:underapproximation}

\begin{figure}[t]
	\centering
	\includegraphics[width=\linewidth,keepaspectratio]{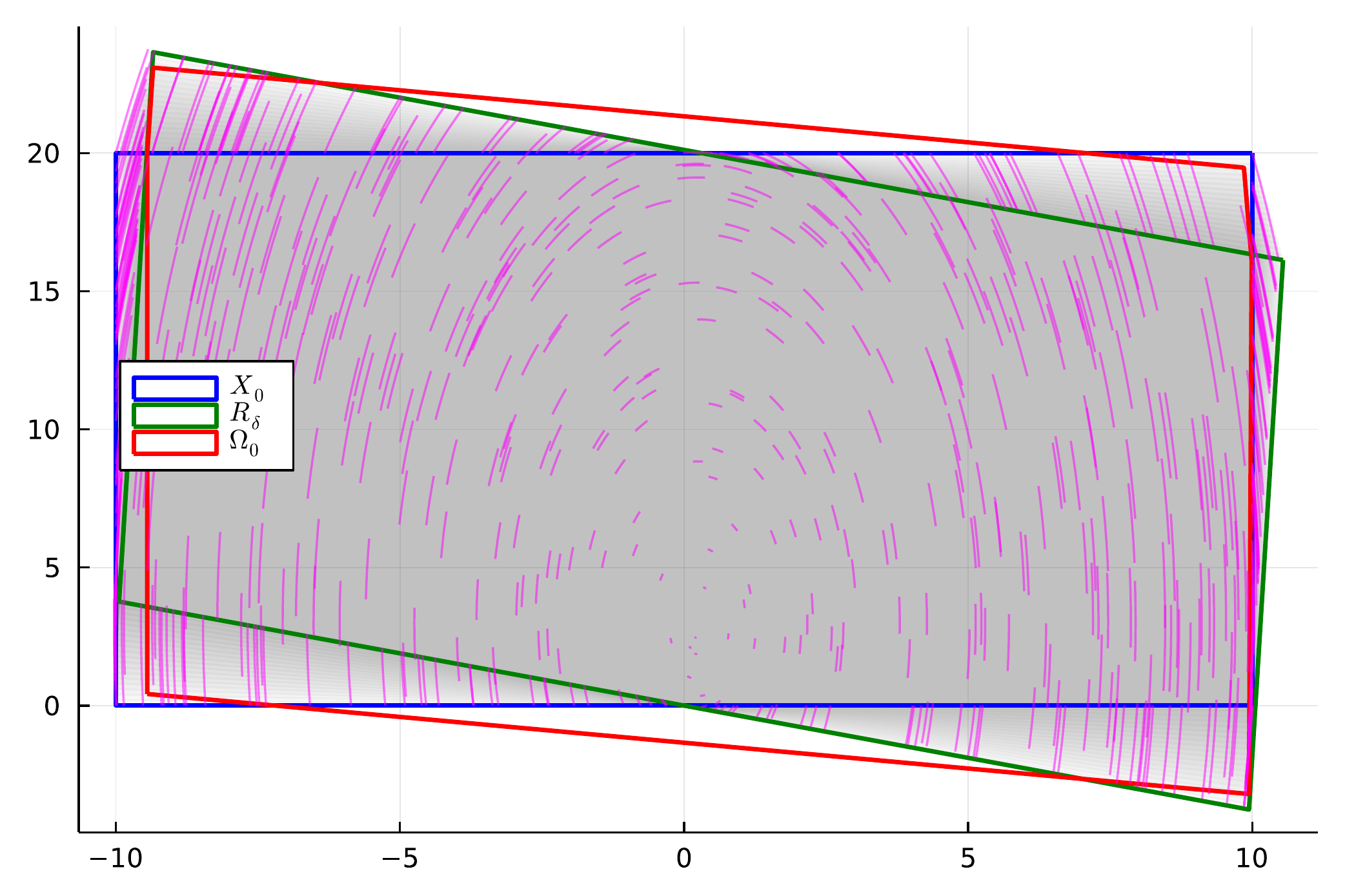}
	\caption{The set $\Omega_0$ (red) obtained with the ``d/dt'' method aimed to underapproximate.
		The system evolves from $\X_0$ (blue) to $\reach(\delta)$ (green) and the reachable set spans the shape of a bow tie, as can be seen from the trajectories (magenta) and some sets $\reach_t$ for uniform $t \in [0, \delta]$ (gray).
		The set $\Omega_0$ contains the unreachable triangles at the top and bottom.}
	\label{fig:underapproximation}
\end{figure}

The authors of the ``d/dt'' method claim that their method can be used to also obtain an underapproximation \cite{AsarinDMB00}.
Recall the equation \eqref{eq:method_homogeneous}.
Instead of bloating the set $\CH(\X_0, \Phi \X_0)$, one can also shrink it.
The authors originally do this by ``inward-pushing'' of the constraints of the polytope, but the idea also naturally extends to a support-function representation of the set.
Denote this shrinking operator with $\ominus$.
The claim is that the set
\begin{equation*}
	\CH(\X_0, \Phi \X_0) \ominus \B_\varepsilon
\end{equation*}
is an underapproximation of the true reachable states.
Unfortunately, the error analysis for $\varepsilon$ only holds for singleton initial states $\X_0 = \{x_0\}$ and is incorrect for proper sets of initial states.
Hence this method does not generally yield a sound underapproximation, which we visualize in Fig.~\ref{fig:underapproximation}.

The problem is that $\varepsilon$ only takes into account the error for the linear interpolation of $x_0$ and $\Phi x_0$, for each $x_0 \in \X_0$; but the convex hull also contains line segments connecting unrelated points $x_0$ and $\Phi x_1$ for $x_0 \neq x_1$.
This becomes apparent if $x_0$ and $x_1$ are sufficiently different, i.e., if $\X_0$ is large, as in the figure.
Take the point $x_0$ as the upper right corner and the point $x_1$ as the upper left corner.
The line segment is the upper edge of the red shape in the figure, which is pushed inward but still falls outside the reachable states.

\section{Homogenization for nondeterministic systems}\label{sec:additional_plots}

\begin{figure}[tb]
	\centering
	\includegraphics[width=0.6\linewidth,keepaspectratio]{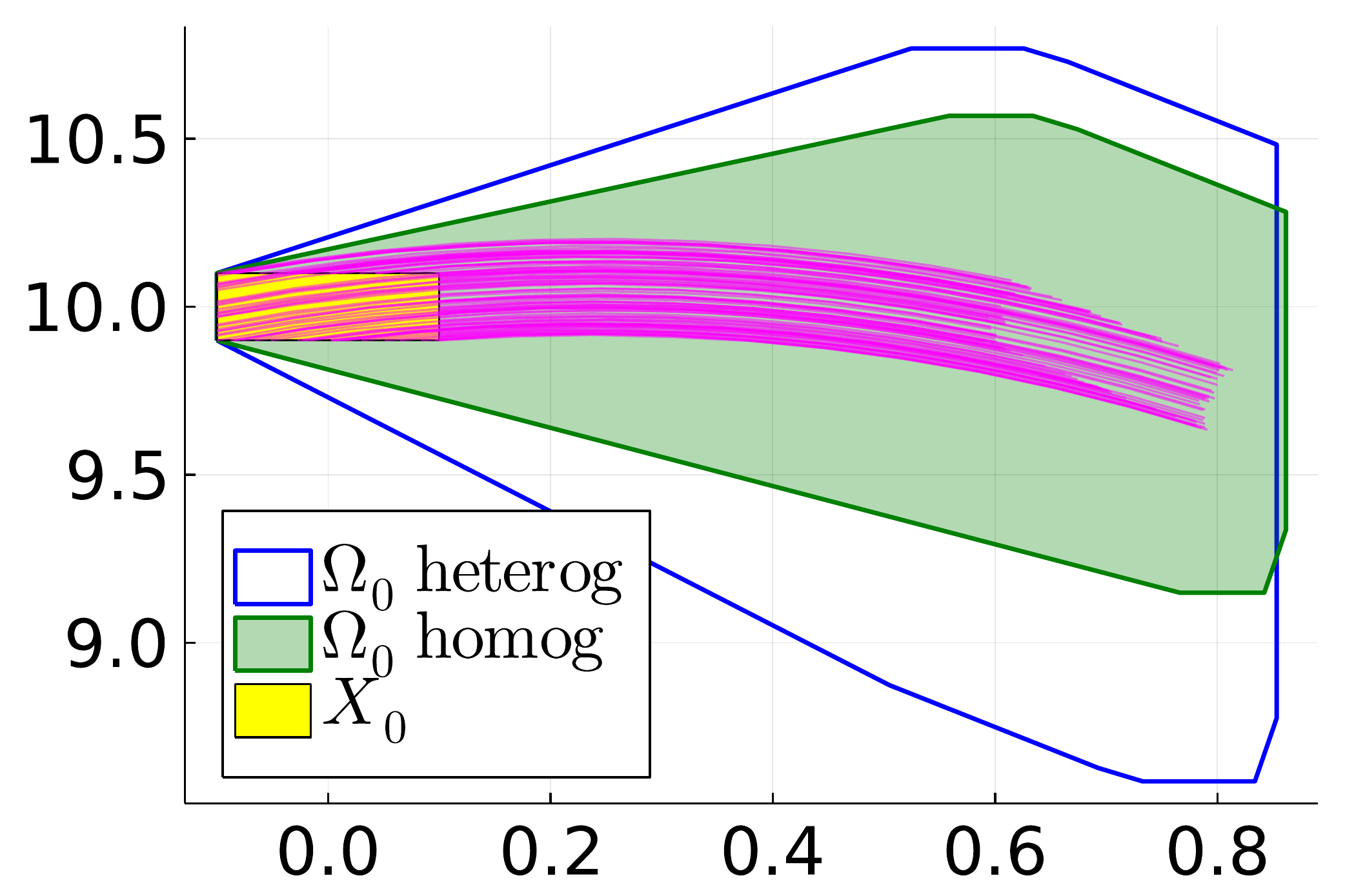}
	\caption{Several trajectories and the sets $\Omega_0$ obtained with the ``forward-only'' method for a nondeterministic heterogeneous system and its partial homogenization, projected to the first two dimensions.}
	\label{fig:homogenize_nondeterministic}
\end{figure}

Fig.~\ref{fig:homogenize_nondeterministic} shows the effect of partially homogenizing a nondeterministic system as described in Section~\ref{sec:homogenize}. It is observed that homogenizing the system improves the precision of the resulting set.

\clearpage

\end{document}